# DIFFUSIVITY IN ONE-DIMENSIONAL GENERALIZED MOTT VARIABLE-RANGE HOPPING MODELS

By P. Caputo and A. Faggionato

*Università Roma Tre and Università "La Sapienza"*

We consider random walks in a random environment which are generalized versions of well-known effective models for Mott variable-range hopping. We study the homogenized diffusion constant of the random walk in the one-dimensional case. We prove various estimates on the low-temperature behavior which confirm and extend previous work by physicists.

**1. Introduction.** Random walks among randomly distributed traps have been proposed as models to study the low-temperature behavior of conductivity in disordered solids in which the Fermi level (set equal to 0 below) lies in a region of strong Anderson localization. In the so-called Mott variable-range hopping model one considers trapping sites $\xi = \{x_i\}$ randomly distributed on $\mathbb{R}^d$, $d \geq 1$, with a given density $\rho$. Each site $x_i$ is marked with a random energy $E_i \in [-1, 1]$, where the variables $E_i$ are independent and identically distributed according to some law $\nu$ on $[-1, 1]$, and are assumed to be independent of $\xi$. The law $\nu$ satisfies $\nu[-E, E] \sim |E|^\delta$ when $E \ll 1$, for some positive constant $\delta$. Then one considers a continuous-time random walk which starts at a given site $x_0$ and jumps from a site $x_i$ to any other site $x_j$ with rate

$$(1.1) \qquad c_{x_i, x_j} = \exp\{-|x_i - x_j| - \beta u(E_i, E_j)\},$$

where $\beta$ is the inverse temperature and the function $u$ is given by

$$(1.2) \qquad u(E_i, E_j) = |E_i| + |E_j| + |E_i - E_j|.$$

The associated random resistor network is obtained by connecting each pair of sites $x_i, x_j$ by the resistor $R_{x_i, x_j} = 1/c_{x_i, x_j}$. At the heuristic level, as predicted by the Einstein relation, the effective conductivity $\sigma$ of the medium

---

can be identified with the diffusion coefficient $D$ associated to the random walk. We refer to [2, 14, 22] and references therein for a justification of the model (1.1) from the physical point of view.

As shown in [2] in the case $\delta = 1$, by means of percolation arguments this model is well suited to explain Mott's law which asserts that the conductivity $\sigma$ should vanish as

$$\log \sigma \sim -\beta^{\delta/(d+\delta)}, \qquad d \geq 2, \tag{1.3}$$

as $\beta \to \infty$. The heuristics behind this behavior can be roughly explained as follows. As $\beta$ increases, the dominant contribution comes from traps $\widetilde{x}_i$ with energy $E_i \in [-\widetilde{E}(\beta), \widetilde{E}(\beta)]$ where $\widetilde{E}(\beta) \to 0$ as $\beta \to \infty$. The traps $\widetilde{\xi} = \{\widetilde{x}_i\}$ are a thinning of the original process $\xi = \{x_i\}$. Since $\nu[-E, E] \sim |E|^{\delta}$, these traps have a density $\rho \nu[-\widetilde{E}(\beta), \widetilde{E}(\beta)] \sim \widetilde{E}(\beta)^{\delta}$ and are typically separated by a distance $\ell(\beta) \sim \widetilde{E}(\beta)^{-\delta/d}$. Percolation ideas valid for $d \geq 2$ allow to argue that these traps contribute roughly $-c_1 \ell(\beta) - c_2 \beta \widetilde{E}(\beta)$ to $\log \sigma$. One obtains the desired estimates by optimizing with the choice $\widetilde{E}(\beta) = \beta^{-d/(\delta+d)}$.

Recently, a mathematical proof of Mott's law for the diffusion coefficient $D$ of the variable-range random walk introduced above was presented in [13, 14] for all dimensions $d \geq 2$.

The situation is quite different in dimension one. As argued in the physics literature [17], here one should have two possible regimes:

$$\text{diffusive regime:} \qquad \log \sigma \sim -\beta, \tag{1.4}$$

$$\text{subdiffusive regime:} \qquad \log \sigma \sim -\infty. \tag{1.5}$$

More precisely, in [17] an upper bound on $\sigma$ is obtained by heuristic arguments when $\delta = 1$ and $\xi = \{x_i\}$ is a Poisson point process of density $\rho$. This upper bound is in agreement with (1.4) when $\rho > 1$ and implies that $\sigma = 0$ when $\rho \leq 1$. As we shall see, the point is that a thinning with arbitrarily small density does not contribute to the conductivity when $d = 1$, that is, nonzero contributions come only from traps with sufficiently high density. We thus have an insulator/conductor transition with critical density $\rho_c = 1$.

The exponential law in (1.4) is in striking contrast with the stretched exponential appearing in the standard Mott law. There is a large literature on the conductivity of one-dimensional or quasi one-dimensional physical systems where the two regimes are compared, see, for example, [1, 18, 19, 21]. Recently, the one-dimensional variable-range hopping model has been also used to study electrical properties of the DNA double helix [23]. All these works focus also on finite-size effects which allow for the onset of a stretched exponential law which is experimentally observed at suitably low temperatures.

Let us now turn to a discussion of our results. We first establish the functional central limit theorem for the diffusively rescaled variable-range random walk. We then characterize its limiting diffusion coefficient $D$ and prove rigorous bounds on $D$. Our bounds will confirm in particular the exponential behavior (1.4) and the transition predicted in [17]. Although the model is spatially one dimensional, the presence of multiple jumps between arbitrarily far traps does not allow for explicit solutions.

Moreover, we will investigate the generalized model which is obtained by replacing the distance $|x_i - x_j|$ in (1.1) by $|x_i - x_j|^\alpha$, with a new parameter $\alpha > 0$. Although the standard case $\alpha = 1$ is the most natural from the physical point of view, we shall see that the case $\alpha < 1$ allows us to produce the various stretched exponential laws discussed in the physics literature. In particular, our main results will establish that if the distribution $\nu$ satisfies $\nu[-E, E] \sim E^\delta$, $\delta > 0$, if $\xi$ is a Poisson point processes with density $\rho$, then one has the following behavior of $\log D$ as a function of $\beta$:

$$(1.6) \qquad \log D \sim \begin{cases} -\beta, & \alpha = 1, \rho > 1, \\ -\beta^{\delta\alpha/(1-\alpha+\delta\alpha)}, & \alpha < 1, \rho > 0, \\ -\infty, & \alpha > 1, \rho > 0 \text{ or } \alpha = 1, \rho \leq 1. \end{cases}$$

We shall observe that for $\alpha \geq 1$ the behavior of the diffusion coefficient $D$ is qualitatively the same as that of the diffusion coefficient $D'$ of the nearest-neighbor random walk on $\xi = \{x_i\}$ with jump rates between consecutive sites $x_i$ and $x_{i+1}$ given by (1.1).

On the other hand, the stretched exponential behavior in the case $\alpha < 1$ will be obtained by optimization arguments partially inspired by the one discussed above for the higher-dimensional case, since in this case one can show that nontrivial contributions come from *any* thinning of the original process. However, we stress that in dimension 1 percolation arguments do not work and in contrast to the higher-dimensional case jumps contributing to $\log \sigma$ can be arbitrarily long if the distances between consecutive points in $\xi = \{x_i\}$ are not bounded from above. In particular, the heuristics given above has to be modified. We will show that the contribution of traps with energy inside the interval $[-\widetilde{E}(\beta), \widetilde{E}(\beta)]$ is no longer $-c_1\ell(\beta) - c_2\beta\widetilde{E}(\beta)$ but rather $-c_1\ell(\beta)^{\alpha/(1-\alpha)} - c_2\beta\widetilde{E}(\beta)$, with $\ell(\beta) \sim \widetilde{E}(\beta)^{-\delta}$. This will give the behavior in (1.6) for $\alpha < 1$, by setting $\widetilde{E}(\beta) = \beta^{-(1-\alpha)/(1-\alpha+\delta\alpha)}$.

Finally, we shall investigate a different but related problem, namely the relaxation speed of the random walk in finite boxes, via spectral gap and isoperimetric estimates. This gives another point of view to discuss the transition from diffusive to subdiffusive behavior. This approach was developed in [6] for the case of dimension $d \geq 2$. Here we describe some finer results that can be obtained in the one-dimensional setting. An important difference with respect to the case $d \geq 2$ is that here we expect that the diffusion coefficient is positive if and only if the finite volume spectral gaps obey diffusive estimates.

1.1. *Model and results.* We consider the following generalization of the variable-range random walk discussed above. Let $\{Z_j, j \in \mathbb{Z}\}$ denote a stationary and ergodic sequence of positive random variables with finite mean. Define the variables $\{x_k, k \in \mathbb{Z}\}$ by

$$(1.7) \qquad x_k = x_{k-1} + Z_{k-1}, \qquad k \in \mathbb{Z}, \; x_0 = 0.$$

Then, we consider the point process $\xi = \{x_k\}$ on $\mathbb{R}$. If $\{Z_j\}$ are i.i.d. then $\xi$ is a renewal process. If $\{Z_j\}$ are i.i.d. exponentially distributed variables then $\xi$ is a homogeneous Poisson process conditioned to contain the origin. We refer to [7] for a basic reference on point processes.

Let $\Theta = \{E_k\}_{k \in \mathbb{Z}}$ denote a family of i.i.d. random variables with values in $[-1, 1]$ with common law $\nu$. We do not make any special assumption on $\nu$. The processes $\xi$ and $\Theta$ are supposed to be independent and we denote by $\omega = (\xi, \Theta)$ the resulting marked point process (we interpret $E_j$ as an energy mark for the point $x_j$). We denote by $\mathbb{P}$ and $\mathbb{E}$ the associated probability measure and expectation.

Given a realization of the environment $\omega$, $X^\omega(t)$ denotes the continuous-time random walk having state space $\xi = \{x_k\}$, starting at the origin and jumping from $x_i$ to $x_j$, $i \neq j$, with jump rate

$$(1.8) \qquad c_{x_i, x_j}(\omega) = \exp\{-|x_i - x_j|^\alpha - \beta u(E_i, E_j)\}.$$

Here $\alpha > 0$ and the function $u$ is assumed to satisfy

$$(1.9) \qquad c_1(|E_i| + |E_j|) \leq u(E_i, E_j) = u(E_j, E_i) \leq c_2(|E_i| + |E_j|)$$

for some positive constants $c_1 \leq c_2$, for any $E_i, E_j \in [-1, 1]$. We set $c_{x_i, x_i}(\omega) = 0$. When $\alpha = 1$ and $u$ is given by (1.2) we are back to the standard model (1.1).

Letting $P^\omega$ denote the law of $X^\omega$, the dynamics is described by the following identities:

$$P^\omega(X^\omega(t + dt) = x_j | X^\omega(t) = x_i) = c_{x_i, x_j}(\omega)\, dt + o(dt), \qquad t \geq 0, \; i \neq j,$$

$$P^\omega(X^\omega(t + dt) = x_i | X^\omega(t) = x_i) = 1 - \sum_{j \,:\, j \neq i} c_{x_i, x_j}(\omega)\, dt + o(dt), \qquad t \geq 0.$$

Equivalently, the random walk $X^\omega$ can be described as follows: after arriving at site $x_i$ the particle waits an exponential time with parameter

$$(1.10) \qquad \lambda_{x_i}(\omega) = \sum_{j \in \mathbb{Z}} c_{x_i, x_j}(\omega)$$

and then jumps to site $x_j$, $j \neq i$, with probability

$$(1.11) \qquad \frac{c_{x_i, x_j}(\omega)}{\lambda_{x_i}(\omega)}.$$

By standard methods (see, e.g., [5] and [14], Appendix A), one can check that the random walk $X^\omega$ is well defined for $\mathbb{P}$-a.a. $\omega$ as soon as $\mathbb{E}(\lambda_0(\omega)) < \infty$. One can easily verify that the above condition is equivalent to requiring

$$(1.12) \qquad \mathbb{E}(\xi[0,1]) < \infty,$$

where $\xi[0,1]$ denotes the number of points of $\xi$ in the interval $[0,1]$. Note that this always holds for $\xi$ a renewal process, since in this case (see Lemma A.1) one has

$$(1.13) \qquad \mathbb{E}(\xi[0,1]^k) < \infty \qquad \forall k \in \mathbb{N}.$$

In order to state our main results we need some further notation. $P^\omega$, the law of $X^\omega$, is a probability measure on the space $D([0,\infty),\xi)$ of right continuous paths having left limits from $[0,\infty)$ to $\xi$, endowed of the Skorokhod topology [4, 10]. We write $E^\omega$ for the associated expectation. We will often identify the random sets $\xi$, $\omega$ with random measures, namely

$$\xi = \sum_{i\in\mathbb{Z}} \delta_{x_i}, \qquad \omega = \sum_{i\in\mathbb{Z}} \delta_{(x_i,E_i)}.$$

Given functions $f = f(\xi)$, $g = g(\omega)$ we set

$$\nabla_y f(\xi) = f(\tau_y \xi) - f(\xi), \qquad \nabla_y g(\omega) = g(\tau_y \omega) - g(\omega), \qquad y \in \xi,$$

where

$$\tau_y \xi = \sum_{i\in\mathbb{Z}} \delta_{x_i - y}, \qquad \tau_y \omega = \sum_{i\in\mathbb{Z}} \delta_{(x_i - y, E_i)}, \qquad y \in \xi.$$

The following functional central limit theorems (FCLTs) hold for all $\alpha > 0$ and for all probability distributions $\nu$ of the energy marks.

THEOREM 1.1. (i) *As $\varepsilon \to 0$, the law $P^\omega_\varepsilon$ of the diffusively rescaled random walk*

$$X^\omega_\varepsilon := (\varepsilon X^\omega(\varepsilon^{-2}t) : t \geq 0),$$

*converges weakly in $\mathbb{P}$-probability to the law of a 1D Brownian motion with diffusion coefficient $D(\beta)$ admitting the variational characterization*

$$(1.14) \qquad D(\beta) = \inf_{g \in L^\infty(\mathbb{P})} \mathbb{E}\left\{ \sum_{i\in\mathbb{Z}} c_{0,x_i}(x_i + \nabla_{x_i} g)^2 \right\}.$$

(ii) *If in addition $\mathbb{E}[(c_{0,x_1})^{-1}] < \infty$, then for $\mathbb{P}$-a.a. environments $\omega$ as $\varepsilon \to 0$ the law $P^\omega_\varepsilon$ converges weakly to the law of the 1D Brownian motion with diffusion coefficient $D(\beta)$.*

We recall that the weak convergence in $\mathbb{P}$-probability stated above simply means that, for any bounded continuous function $F$ on the path space $D([0, \infty), \mathbb{R})$ endowed with the Skorokhod topology, the random variable

$$\omega \rightarrow E_\varepsilon^\omega (F(X_\varepsilon^\omega))$$

with $E_\varepsilon^\omega$ denoting expectation with respect to $P_\varepsilon^\omega$, converges in $\mathbb{P}$-probability to the expectation of $F(W)$, where $W$ is a 1D Brownian motion with diffusion coefficient $D(\beta)$, that is, the variance of $W_t$ equals $D(\beta)t$.

The first part of the above theorem will be obtained along the lines of classical homogenization results [8, 16]. It is derived also in [14] under an additional finite moment condition. The second part of the theorem is an almost sure FCLT and its proof will be based on the construction of the so-called corrector field.

We point out that the diffusion coefficient $D(\beta)$ characterized by (1.14) could be zero, thus implying the subdiffusive behavior of the random walk. If not zero, it is relevant to analyze the behavior of $D(\beta)$ as $\beta$ goes to infinity in the same spirit of Mott's law. In the next theorem, we answer these issues by giving upper and lower bounds on $D(\beta)$ under suitable assumptions. Before stating our results, it is convenient to adopt the following notation: given functions $f(x)$ and $g(x)$, $x$ positive number, we write

$$f(x) \sim g(x), \qquad x \searrow 0 [x \nearrow \infty],$$

if there exist $c_1, c_2 > 0$ such that

$$c_1 g(x) \leq f(x) \leq c_2 g(x)$$

for all $x$ in a neighborhood of 0 $[\infty]$.

THEOREM 1.2. *Suppose that $Z_j$ has finite second moment, that is, $\mathbb{E}[Z_j^2] < \infty$.*

(1) *If $\mathbb{E}[\exp{(Z_j^\alpha)}] < \infty$, then there exist $C, \kappa > 0$ such that for all $\beta \geq 0$*

$$(1.15) \qquad\qquad D(\beta) \geq C \exp{[-\kappa\beta]}.$$

(2) *Suppose that $\alpha \geq 1$ and that $\{Z_j\}$ are positive i.i.d. random variables. Then*

$$(1.16) \qquad\qquad D(\beta) > 0 \quad \Longleftrightarrow \quad \mathbb{E}[\exp(Z_j^{\,\alpha})] < \infty.$$

*If in addition to the above bound $\nu$ has no mass at zero, that is, $\nu[-E, E] \rightarrow 0$ as $E \searrow 0$, then*

$$(1.17) \qquad\qquad C_1 \exp[-\kappa_1\beta] \leq D(\beta) \leq C_2 \exp[-\kappa_2\beta]$$

*for suitable positive constants $C_1, C_2, \kappa_1, \kappa_2$, for all $\beta \geq 0$.*

(3) *Suppose that $\alpha < 1$ and that $\{Z_j\}$ are positive i.i.d. random variables.*

(i) If $\mathbb{E}[\exp(\gamma Z_j{}^\alpha)] = \infty$ for some $\gamma \in (0,1)$, then $D(\beta) = 0$.

(ii) If $\mathbb{E}[\exp(\varepsilon Z_j)] < \infty$ for some $\varepsilon > 0$ and

$$(1.18) \qquad \nu[-E, E] \sim E^\delta, \qquad E \searrow 0,$$

for some constant $\delta > 0$, then

$$(1.19) \quad C_1 \exp[-\kappa_1 \beta^{\delta\alpha/(1-\alpha+\delta\alpha)}] \leq D(\beta) \leq C_2 \exp[-\kappa_2 \beta^{\delta\alpha/(1-\alpha+\delta\alpha)}]$$

for suitable positive constants $C_1, C_2, \kappa_1, \kappa_2$, for all $\beta \geq 0$.

The upper bounds in (1.17) and (1.19) will be obtained by choosing suitable test functions in the variational characterization (1.14). The proof of the other bounds on $D(\beta)$ will be based on comparison with suitable nearest-neighbor random walks for whose diffusion coefficients one can derive explicit expressions. For instance, to define a nearest-neighbor walk on $\xi$ with a diffusion coefficient smaller than $D(\beta)$ we shall simply cut all jumps of the form $x_i \to x_j$ with $|i - j| \neq 1$. This produces the random walk $Y^\omega$ which will be used in the proof of (1.15). On the other hand, the construction of a nearest-neighbor walk on $\xi$ whose diffusion coefficient is larger than $D(\beta)$ is more tricky. This will be needed for the proof of points (2) and (3)(i) in Theorem 1.2. Here the electric network representation and the associated monotonicity laws [9] will turn out to be a very useful guide for our comparison arguments. A convenient construction in the electric network context has been outlined in [17].

It will be clear from our proof that the diffusion coefficient of the nearest-neighbor walk $Y^\omega$ defined above is positive iff the condition in the right-hand side of (1.16) is satisfied and, under this assumption it satisfies the same bounds appearing in (1.17). Hence, we can conclude that for $\alpha \geq 1$ and at least for renewal point processes $\xi$ the Mott variable-range hopping is essentially a nearest-neighbor hopping. If $\alpha < 1$ on the other hand, the proof of the stretched exponential estimates in point (3)(ii) above will make clear that the main contribution to diffusion comes from jumps between sites $x_i$ with small energy marks, which are typically not nearest neighbors.

Let us now draw some consequences of Theorem 1.2 in specific cases. Various other cases can be considered in the way described below.

1.1.1. *Poissonian case.* This is the case where $\{Z_j\}$ are i.i.d. exponentials of parameter $\lambda > 0$. Here we can easily collect the estimates in Theorem 1.2 to obtain the picture anticipated in (1.6), where $\rho = \lambda$. In particular, the critical point for diffusivity appears at $\alpha = \alpha_c = 1$, $\lambda = \lambda_c = 1$ and the system is subdiffusive at the critical point.

1.1.2. *Diluted lattice.* Suppose the $\{Z_j\}$ are i.i.d. geometric random variables with parameter $p \in (0, 1)$, that is, $\mathbb{P}(Z_j = k) = p(1-p)^{k-1}, k \in \mathbb{N}$. This corresponds to the point process $\xi$ obtained by deleting independently each point of $\mathbb{Z} \setminus \{0\}$ with probability $1 - p$. It is easily seen that here one has the same picture as in the Poissonian case described above with the correspondence $p = 1 - e^{-\lambda}$.

1.1.3. *Weibull distribution.* Here we take $\{Z_j\}$ to be i.i.d. with a Weibull distribution, that is, $\mathbb{P}(Z_j > t) = e^{-\lambda t^\tau}, t \geq 0$, for some $\lambda > 0, \tau > 0$. This generalizes the Poissonian case ($\tau = 1$). It follows from Theorem 1.2 that the critical point for diffusivity is given by $\alpha = \alpha_c = \tau$, and $\lambda = \lambda_c = 1$. Thanks to point (2) in the theorem we know that the system is subdiffusive at the critical point at least in the case $\tau \geq 1$. For $\tau < 1$ we cannot exclude positive diffusion constant at the critical point.

1.1.4. *Gaussian distribution.* We may take $\{Z_j\}$ to be i.i.d. with $Z_j$ distributed as $|Y|$ where $Y$ is a normal random variable with mean 0 and variance $\sigma^2$. Since, apart from polynomial corrections, $\mathbb{P}(Y \geq t)$ behaves as $e^{-t^2/(2\sigma^2)}$ for $t$ large, it is simple to check that in this case the asymptotic behavior of $D(\beta)$ corresponds to the case of Weibull distribution with $\tau = 2$, $\lambda = 1/(2\sigma^2)$.

1.1.5. *Random variables with polynomial tails.* If we take $\{Z_j\}$ to be i.i.d. with $\mathbb{P}(Z_j > t) \sim t^{-r} \; \forall t \geq 1$, for some $r > 2$ (we require $r > 2$ so that $\mathbb{E}[Z_j^2] < \infty$), then it is immediate to use points (2) and (3)(i) in the theorem to see that the system is subdiffusive for all $\alpha > 0$ and all $r > 2$.

1.2. *Diffusivity via spectral gap and isoperimetric bounds.* As in [6] our finite volume estimates will be related to the geometry of the point process only and the temperature plays essentially no role. Therefore, without loss of generality, we set $\beta = 0$ for this subsection. We stress that the results presented here are proven for almost all environments under the assumption that $\xi$ is a renewal process. As a convention, whenever we state that some event $E_L$ involving the parameter $L$ holds $\mathbb{P}$-a.s. we mean that with $\mathbb{P}$-probability 1 there exists $L_0 = L_0(\xi) < \infty$ such that the event $E_L$ occurs for all $L \in \mathbb{N}$, with $L \geq L_0$.

Let $\xi = \{x_k, k \in \mathbb{Z}\}$ denote as usual the renewal process generated by (1.7). We write, for every $L > 0$, $\xi_L := \xi \cap \Lambda_L$, where $\Lambda_L = [-L/2, L/2]$, for the process inside the segment of length $L$.

We recall the definition of Cheeger's isoperimetric constant $\Phi_L(\xi)$:

$$(1.20) \qquad \Phi_L(\xi) := \min_{U \subset \xi_L \,:\, \#(U) \leq 1/2 \#(\xi_L)} I_U(\xi),$$

$$(1.21) \qquad I_U(\xi) := \frac{1}{\#(U)} \sum_{x \in U, y \in \xi_L \setminus U} e^{-|x-y|^\alpha},$$

where $\#(U)$ stands for the number of points in $U$.

In order to estimate the isoperimetric constant $\Phi_L(\xi)$ it is convenient to introduce the random variable $\zeta_L$ defined as the maximal distance between consecutive points in $\xi_L$:

$$(1.22) \qquad \zeta_L = \max\{x_{k+1} - x_k : x_k, x_{k+1} \in \xi_L\}.$$

Then the behavior of $\Phi_L$ is mostly determined by $\zeta_L$.

PROPOSITION 1.3. *For any $\gamma \in (0,1)$, there exists a positive constant $C$ such that*

$$(1.23) \qquad \frac{1}{CL}e^{-\zeta_L^\alpha} \le \Phi_L(\xi) \le \frac{C(\log L)^2}{L}e^{-d(\alpha)\zeta_{\gamma L}^\alpha}, \qquad \mathbb{P}\text{-}a.s.,$$

*where $d(\alpha) := 1 \wedge 3^{\alpha-1}$. If there exists $\varepsilon > 0$ such that $\mathbb{P}(Z_i > \varepsilon) = 1$ then (1.23) holds without the $(\log L)^2$ factor in the right-hand side.*

Since the $Z_i$ are independent, the a.s. behavior of $\zeta_L$ can be characterized in terms of the tail

$$\psi(t) = \mathbb{P}(Z_i > t).$$

We note that $\lim_{L\uparrow\infty} \zeta_L = \infty$ a.s. if $Z_i$'s are unbounded and $\lim_{L\uparrow\infty} \zeta_L < \infty$ a.s. if $Z_i$'s are bounded [11], Chapter 3. Hence, from Proposition 1.3, $\Phi_L(\xi) \ll L^{-1}$ a.s. if $Z_i$ are unbounded, while $\Phi_L \sim L^{-1}$ a.s. if $Z_i$'s are bounded and bounded away from zero. In [6] the bound $\Phi_L \sim L^{-1}$ was shown to be characteristic of the diffusive regime for the spectral gap of the random walk in dimension $d \ge 2$. When $d = 1$, we shall see that Proposition 1.3 and the results of the previous section imply that this is no longer true. Namely, the spectral gap can have diffusive behavior while $\Phi_L \ll L^{-1}$.

Proposition 1.3 and standard results concerning the a.s. asymptotic behavior of maxima of i.i.d. random variables (see [15], Section 4, [11], Section 3.5) allow to derive rather fine estimates on $\Phi_L$. In Section 6 we discuss a general method and, as an application, we derive the following criterion which covers several interesting cases discussed below.

THEOREM 1.4. *Let $\{Z_i, i \in \mathbb{Z}\}$ be i.i.d. positive variables with finite mean $\mu$. Let $a, b : \mathbb{R}_+ \to \mathbb{R}_+$ be two nonnegative, nondecreasing functions satisfying $\psi(a(n)), \psi(b(n)) \to 0$, and $n\psi(b(n)) \to \infty$ as $n \to \infty$ and such that*

$$(1.24) \qquad \sum_{n\in\mathbb{N}} \psi(a(n)) < \infty,$$

$$(1.25) \qquad \sum_{n\in\mathbb{N}} \psi(b(n))e^{-n\psi(b(n))} < \infty.$$

*Then, for any $\gamma \in (0, 1)$, there exists a positive constant $C = C(\gamma) < \infty$ such that*

$$(1.26) \quad \frac{1}{C} \frac{e^{-a(L/(\gamma\mu))^\alpha}}{L} \leq \Phi_L \leq C(\log L)^2 \frac{e^{-d(\alpha)b(\gamma L/\mu)^\alpha}}{L}, \qquad \mathbb{P}\text{-}a.s.,$$

*where $d(\alpha) = 1$ if $\alpha \geq 1$ and $d(\alpha) = 3^{\alpha-1}$ if $\alpha < 1$. Moreover, if the $Z_i$'s are bounded one can take $a$ and $b$ constant in (1.26). Finally, if the $Z_i$'s are bounded away from zero then the estimate in the right-hand side of (1.26) holds without the logarithmic correction.*

Let us now look at the Poincaré or spectral gap inequality associated to the random walk in $\Lambda_L$. The Poincaré constant $\gamma(L) = \gamma(L, \xi)$ is defined by

$$(1.27) \qquad \gamma(L) := \sup_{f: \xi_L \to \mathbb{R}} \frac{1}{\#(\xi_L)} \frac{\sum_{x,y \in \xi_L} (f(x) - f(y))^2}{\sum_{x,y \in \xi_L} e^{-|x-y|^\alpha} (f(x) - f(y))^2}.$$

The spectral gap is given by $\text{gap}(L) := \gamma(L)^{-1}$. We recall that $\text{gap}(L)$ coincides with the smallest nonzero eigenvalue of the nonnegative matrix $-\mathcal{L}_L$, where $\mathcal{L}_L$ is the Markov generator of the variable-range random walk confined to the segment $\Lambda_L$, defined by its action on vectors $f: \xi_L \to \mathbb{R}$ as follows:

$$(1.28) \qquad \mathcal{L}_L f(x) = \sum_{y \in \xi_L} e^{-|x-y|^\alpha} (f(y) - f(x)), \qquad x \in \xi_L.$$

In the diffusive regime we expect the Laplacian-like behavior $\text{gap}(L) \sim L^{-2}$. In view of our previous discussion this cannot be obtained in general through Cheeger's inequality asserting that $\gamma(L) \leq 8\Phi_L^{-2}$; cf. [6]. However, a simple alternative is available in dimension one and we are able to prove that $\text{gap}(L) \geq cL^{-2}$, $\mathbb{P}$-a.s. as soon as $\mathbb{E}[e^{Z_1^\alpha}] < \infty$. On the other hand, the simple estimate $\gamma(L) \geq \frac{1}{2}\Phi_L^{-1}$ [which follows from (1.27) by restricting $f$ to be an indicator function] allows to obtain upper bounds on the spectral gap from Theorem 1.4. We can summarize these facts in the following:

COROLLARY 1.5. *Let $\{Z_i, i \in \mathbb{Z}\}$ be i.i.d. positive variables with finite mean $\mu$.*

1. *If $\mathbb{E}[e^{Z_1^\alpha}] < \infty$ then there exists $c > 0$ such that*

$$(1.29) \qquad \text{gap}(L) \sim L^{-2}, \qquad \mathbb{P}\text{-}a.s.$$

2. *Let the function $b: \mathbb{R}_+ \to \mathbb{R}_+$ be as in Theorem 1.4. Then, for any $\gamma \in (0, 1)$, there exists a positive constant $C = C(\gamma) < \infty$ such that*

$$(1.30) \qquad \text{gap}(L) \leq C(\log L)^2 \frac{e^{-d(\alpha)b(\gamma L/\mu)^\alpha}}{L}, \qquad \mathbb{P}\text{-}a.s.,$$

*where $d(\alpha) = 1 \wedge 3^{\alpha-1}$. If $\mathbb{P}(Z_i \leq \delta) = 0$ for some $\delta > 0$ the last estimate holds without the logarithmic correction.*

1.2.1. *Examples.* Let us conclude this discussion with some details on specific examples. We start with the *Poisson case*, that is, $\{Z_j\}$ are i.i.d. exponential variables of parameter $\lambda > 0$.

We can apply Theorem 1.4 with the following choice of the functions $a, b$:

$$(1.31) \qquad a(t) = \frac{1}{\lambda} \log[t(\log t)^p], \qquad b(t) = \frac{1}{\lambda} \log[t/(p \log \log t)]$$

for some $p > 1$ and all $t$ large enough. Here $\psi(t) = e^{-\lambda t}$ and one can easily check that both (1.24) and (1.25) are satisfied, since with this choice we have

$$(1.32) \qquad \psi(a(n)) = \frac{1}{n(\log n)^p}, \qquad \psi(b(n)) = \frac{p \log \log n}{n},$$

$$\psi(b(n)) e^{-n\psi(b(n))} = \frac{p \log \log n}{n(\ln n)^p}.$$

If $\alpha = 1$, from (1.26) we then obtain in particular that for every $\lambda > 0$ there exists some constant $c > 0$ such that

$$(1.33) \qquad (\log L)^{-c} L^{-1-1/\lambda} \leq \Phi_L \leq (\log L)^c L^{-1-1/\lambda}, \qquad \mathbb{P}\text{-a.s.}$$

Since $\Phi_L$ behaves like $L^{-1-1/\lambda}$ apart from logarithmic corrections, Corollary 1.5 implies in particular that $\text{gap}(L) \sim L^{-2}$, $\mathbb{P}$-a.s. whenever $\lambda > 1$, while $\text{gap}(L) \leq (\log L)^c L^{-1-1/\lambda}$, when $\lambda < 1$. Hence, we have a transition at $\lambda = 1$ from diffusive to subdiffusive relaxation.

If $\alpha > 1$ then we easily see from Theorem 1.4 and Corollary 1.5 that both $\text{gap}(L)$ and $\Phi_L$ decay to zero almost surely, faster than any inverse power of $L$.

When $\alpha \in (0, 1)$, we know by Corollary 1.5 that $\text{gap}(L) \sim L^{-2}$, $\mathbb{P}$-a.s. while $L \Phi_L \to 0$ as $L \to \infty$ since $Z_j$'s are unbounded.

The same arguments apply to the *diluted lattice*, that is, to the case of geometric random variables with parameter $p$, with the correspondence $p = 1 - e^{-\lambda}$. Moreover, the extension to *Weibull distributions* is also straightforward. Here we have $\psi(t) = e^{-\lambda t^\tau}$, for some $\lambda > 0, \tau > 0$. We simply replace the functions in (1.31) by $a(t) \to a(t)^{1/\tau}$ and $b(t) \to b(t)^{1/\tau}$. Then it follows that (1.24) and (1.25) are satisfied, just as in the case $\tau = 1$ discussed above. In conclusion, thanks to Corollary 1.5, one has that $\text{gap}(L)$ shrinks faster than any inverse power of $L$ for $\alpha > \tau$, while $\text{gap}(L) \sim L^{-2}$ for $\alpha < \tau$. At $\alpha = \tau \geq 1$ one finds again a transition at $\lambda = 1$. Indeed, a direct inspection of the bounds in (1.26) shows that for $\alpha = \tau \geq 1$ one has again the estimate (1.33) for every $\lambda > 0$. In particular, Corollary 1.5 shows a transition from $\text{gap}(L) \sim L^{-2}$ (when $\lambda > 1$) to $\text{gap}(L) \leq L^{-2-\varepsilon}$ for some $\varepsilon > 0$ (when $\lambda < 1$). When $\alpha = \tau < 1$ we believe that the same should hold. However, since now $d(\alpha) = 3^{\alpha-1} < 1$, then the right-hand side of (1.33) has to be modified with replacement of $(\log L)^c L^{-1-1/\lambda}$ by $(\log L)^c L^{-1-3^{\alpha-1}/\lambda}$. Therefore, by Corollary 1.5, if $\alpha = \tau < 1$, $\mathrm{gap}(L) \ll L^{-2}$ for $\lambda < 3^{\alpha-1}$ while $\mathrm{gap}(L) \sim L^{-2}$ for $\lambda > 1$, thus we are only able to say that a transition occurs for $\lambda$ somewhere between $3^{\alpha-1}$ and 1. Note that, when $\tau = 2$, the same behavior is produced in the *Gaussian case*, that is, by random variables $Z_i = |Y_i|$, where $Y_i$ are normally distributed with mean 0 and variance $\sigma^2 = 1/(2\lambda)$.

Finally, one may consider *variables with a polynomial tail*, such as, for example, $\psi(t) = t^{-r}$, with some $r > 1$, for all $t \geq 1$. Theorem 1.4 allows to obtain detailed estimates on $\Phi_L$ as follows. Fix $p > 1$ and define, for $t$ large enough,

$$a(t) = t^{1/r} (\log t)^{p/r}, \qquad b(t) = \left( \frac{t}{p \log \log t} \right)^{1/r}.$$

The above sequences are increasing for $t$ large enough, moreover one can check that (1.32) holds. In particular, Theorem 1.4 can be used to produce the estimates, for any $\varepsilon > 0$

$$\exp[-L^{\alpha/r} (\log L)^{\alpha/r+\varepsilon}] \leq \Phi_L \leq \exp[-L^{\alpha/r} (\log L)^{-\varepsilon}], \qquad \mathbb{P}\text{-a.s.}$$

1.3. *Overview of the following sections.* The paper is organized as follows. The proof of Theorem 1.1 is given in Section 2. The proof of Theorem 1.2 requires several distinct tools which are given separately in Sections 3, 4 and 5. In Section 3 we introduce nearest-neighbor random walks, characterize their diffusion constants and prove the lower bound (1.15) and the lower bound in the left-hand side of (1.19). In Section 4 we prove the equivalence (1.16) and statement (3)(i) of Theorem 1.2. In Section 5 we prove the upper bounds in (1.17) and (1.19). Finally, in Section 6 we prove our results on isoperimetric constants and spectral gaps.

**2. Proof of Theorem 1.1.** In order to prove Theorem 1.1, it is convenient to fix some notation. We denote by $\mathcal{N}_0$ the space of doubly–infinite sequences $\omega = \{x_i, E_i\}_{i \in \mathbb{Z}}$, where $x_i < x_{i+1}$, $x_0 = 0$ and $E_i \in [-1, 1]$. The topology of $\mathcal{N}_0$ is the one making the map

$$\mathcal{N}_0 \ni \omega \to \{(x_{i+1} - x_i, E_i)\}_{i \in \mathbb{Z}} \in ((0, \infty) \times [-1, 1])^{\mathbb{Z}}$$

a topological homeomorphism. We recall that, given $\omega$ as above, $\xi$ denotes its spatial projection, namely $\xi = \{x_i\}_{i \in \mathbb{Z}}$.

The annealed FCLT given in part (i) can be derived by applying Theorem 2.1 of [8] and then using a time–change argument. We only sketch the proof, since these methods are rather standard. First, we consider the discrete-time random walk $\{\hat{X}^\omega(n)\}$ starting at the origin and jumping from $x_i$ to $x_j$, $j \neq i$, with probability $c_{x_i,x_j}(\omega)/\lambda_{x_i}(\omega)$. The associated *environment viewed from the particle* $\omega_n := \tau_{\hat{X}^\omega(n)} \omega$ is reversible (see Lemma 1(i) in [14]) and ergodic with respect to the probability measure $\mathbb{Q}$ defined as $d\mathbb{Q}(\omega) = (\lambda_0(\omega)/\mathbb{E}(\lambda_0))d\mathbb{P}(\omega)$. Then, we note that the function $f(\omega, \omega')$ defined as $x$ if $\omega' = \tau_x\omega$ and 0 otherwise is $\mathbb{P} \times \mathbb{P}$-a.s. (and therefore $\mathbb{Q} \times \mathbb{Q}$-a.s.) well defined on $\mathcal{N}_0 \times \mathcal{N}_0$ if

$$(2.1) \qquad \tau_x\omega \neq \omega \qquad \forall x \in \xi, \ \mathbb{P}\text{-a.s.}$$

The above condition is violated only if $\xi = a\mathbb{Z}$ and $E_i \equiv E$ for suitable constants $a, E$. Since in this case the analysis of the random walk $X^\omega$ becomes trivial, without loss of generality we can assume (2.1). Then $f$ is a well-defined and antisymmetric function. Moreover, due to our basic assumption (1.12),

$$
\begin{aligned}
(2.2) \qquad E(f(\omega_0, \omega_1)^2) &= \int \mathbb{Q}(d\omega) \sum_{i \in \mathbb{Z}} \frac{c_{0,x_i}(\omega)}{\lambda_0(\omega)} x_i^2 \\
&= \mathbb{E}(\lambda_0)^{-1} \int \mathbb{P}(d\omega) \sum_{i \in \mathbb{Z}} c_{0,x_i}(\omega) x_i^2 < \infty,
\end{aligned}
$$

where the first member denotes the expectation with respect to the process $\{\omega_n\}$, $\omega_0$ being chosen with distribution $\mathbb{Q}$. In particular, all the conditions of Theorem 2.1 in [8] are satisfied. Since $\hat{X}^\omega(n) = \sum_{j=0}^{n-1} f(\omega_j, \omega_{j+1})$, this theorem implies that the diffusively rescaled random walk $\hat{X}_\varepsilon^\omega := (\varepsilon\hat{X}^\omega([\varepsilon^{-2}t]) : t \geq 0)$ converges weakly in $\mathbb{Q}$-probability to the 1D Brownian motion with diffusion coefficient $\hat{D}(\beta)$ defined as (see [8], (2.28))

$$(2.3) \qquad \hat{D}(\beta) = E(f(\omega_0, \omega_1)^2) - 2\|\varphi\|_{-1, L^2(\mathbb{Q})}^2,$$

where $\varphi(\omega) = \sum_{i \in \mathbb{Z}} (c_{0,x_i}(\omega)/\lambda_0(\omega)) x_i$ and $\|\varphi\|_{-1, L^2(\mathbb{Q})}$ denotes the $H_{-1}$-norm of $\varphi$ with respect to the symmetric operator $Lg(\omega) = \sum_{i \in \mathbb{Z}} (c_{0,x_i}(\omega)/\lambda_0(\omega)) \times (g(\tau_i\omega) - g(\omega))$ defined on $L^2(\mathbb{Q})$. Since the family $B_\infty$ of bounded Borel functions is dense in $L^2(\mathbb{Q})$, (2.3) reads

$$\hat{D}(\beta) = E(f(\omega_0, \omega_1)^2) - 2 \sup_{g \in B_\infty} (2\langle \varphi, g \rangle_{L^2(\mathbb{Q})} - \langle g, -Lg \rangle_{L^2(\mathbb{Q})}).$$

By standard algebra together with Lemma 1(i) in [14], it is simple to prove that the above right-hand side coincides with $\mathbb{E}(\lambda_0)^{-1}D(\beta)$, with $D(\beta)$ defined as in (1.14). At this point, part (i) of Theorem 1.1 can be obtained by thinking of $X^\omega$ as obtained from $\hat{X}^\omega$ by a time change (see, e.g., [14], page 13) and then applying Theorem 17.1 in [4] as done in the proof of Theorem 4.5 in [8].

Let us now prove the almost sure FCLT stated in Theorem 1.1(ii). Again, we first prove that the diffusively rescaled random walk $\hat{X}_\varepsilon^\omega$ converges to the 1D Brownian motion with diffusion coefficient $\hat{D}(\beta)$ and then extend the result to $X^\omega$ by means of time–change arguments. The proof for $\hat{X}^\omega$ is based on the construction of the so-called corrector field $\chi(\omega, x_i)$, a random variable at each point $x_i$ such that $\xi \ni x \to \varphi(\omega, x) := x + \chi(\omega, x)$ is harmonic with respect to the random walk $\hat{X}^\omega$. This approach is by now rather standard and we shall only point out the main steps required together with the relevant literature; see, for example, [3, 20] for recent accounts. As in these works the FCLT is derived from a sublinearity property of the corrector together with an application of the Lindeberg–Feller CLT for martingales. It is important to note that the 1D setting allows for significant simplification. On the other hand, the long-range nature and the intrinsically nondeterministic structure of the state space of our walk do require some modifications of the standard arguments.

In order to construct the corrector, we introduce the Hilbert space $\mathcal{H}$ of *square-integrable forms* vanishing at the origin, that is, of measurable functions $u : \mathcal{N}_0 \times \mathbb{R} \to \mathbb{R}$ such that $u(\omega, 0) = 0$, $\mathbb{P}$-a.s. and

$$(2.4) \qquad \|u\|_2 := \sum_{i \in \mathbb{Z}} \mathbb{E}[c_{0,x_i}(\omega) u(\omega, x_i)^2] < \infty.$$

Note that, from our basic assumption (1.12), we know that $u(\omega, x) = x$ is an element of $\mathcal{H}$. Let $\mathcal{H}_\Delta \subset \mathcal{H}$ denote the closure of the linear subspace of gradient functions. Namely, $\psi \in \mathcal{H}_\Delta$ iff there exists a sequence of bounded measurable functions $g_n : \mathcal{N}_0 \to \mathbb{R}$ such that $\|\psi - \nabla g_n\|_2 \to 0$, $n \to \infty$, with $\nabla g_n(\omega, x)$ defined as $g_n(\tau_x \omega) - g(\xi)$, for $x \in \xi$, and $0$ otherwise. Elements of $\mathcal{H}_\Delta$ are called *potential forms*. Its orthogonal complement $\mathcal{H}_\nabla^\perp$ is the space of *solenoidal forms*. It can be checked that the standard theory applies. Namely:

(1) Potential forms $u \in \mathcal{H}_\Delta$ are curl-free, that is, they satisfy the cocycle property: for $\mathbb{P}$-a.a. $\omega$, for every $x \in \xi$, for any finite set of points $y_0, \dots, y_m \in \xi$ with $y_0 = y_m$ one has

$$(2.5) \qquad \sum_{j=0}^{m-1} u(\tau_{y_j} \omega, y_{j+1} - y_j) = 0.$$

[The above property trivially holds if $u = \nabla g(\omega, x)$ and extends by continuity to all $\mathcal{H}_\Delta$.]

(2) A form $u \in \mathcal{H}$ is solenoidal if and only if it is divergence free, that is, for $\mathbb{P}$-a.a. $\omega$ it holds

$$(2.6) \qquad \operatorname{div} u(\omega) := \sum_{x \in \xi} c_{0,x}(\omega)(u(\omega, x) - u(\tau_x \omega, -x)) = 0.$$

The above characterization is implied by the identity:

$$\mathbb{E}\left(\sum_{x \in \xi} c_{0,x}(\omega) \nabla g(\omega, x) u(\omega, x)\right) = -\mathbb{E}\left(\sum_{x \in \xi} c_{0,x}(\omega) g(\omega, x) \operatorname{div} u(\omega)\right),$$

where $g$ is a bounded measurable function on $\mathcal{N}_0 \times \mathbb{R}$ and $u \in \mathcal{H}$.

Finally, (3) a form $\varphi$ which is both curl-free and divergence free must be harmonic, that is, for $\mathbb{P}$-a.a. $\omega$, for any $x \in \xi$ we must have

$$(2.7) \qquad \varphi(\omega, x) = \sum_{y \in \xi} c_{x,y}(\omega)(\varphi(\omega, y) - \varphi(\omega, x)).$$

Indeed, due to (2.5) it must be $\varphi(\omega, z) + \varphi(\tau_z \omega, -z) = 0$ for $\mathbb{P}$-a.a. $\omega$, for any $z \in \xi$. This together with (2.6) implies (2.7) for $x = 0$. The general case follows from the covariance of jump rates with respect to translations.

We can finally define the corrector field $\chi$. Consider the form $u$ defined by $u(\omega, x) = x$. Let $\pi : \mathcal{H} \to \mathcal{H}_\nabla$ be the orthogonal projection on potential forms and define the corrector field by $\chi := \pi(-u)$. We see that the form $\varphi = u + \chi$ is curl-free ($u$ is clearly curl-free and $\chi$ is potential). Moreover, by construction $\varphi$ is solenoidal and therefore div $\varphi = 0$. It follows that $\varphi$ is harmonic as in (2.7). We can list the following properties of the corrector $\chi$:

LEMMA 2.1. *(1)* $\chi \in \mathcal{H}$, *that is,* $\chi(\omega, 0) = 0$ *for* $\mathbb{P}$-*a.a.* $\omega$ *and* $\mathbb{E}[\sum_{i \in \mathbb{Z}} c_{0,x_i}(\omega) \chi(\omega, x_i)^2] < \infty$.

*(2) For any* $i \in \mathbb{Z}$, *the map* $\mathcal{N}_0 \ni \omega \to \chi(\omega, x_i)$ *is in* $L^1(\mathbb{P})$ *and* $\mathbb{E}[\chi(\omega, x_i)] = 0$.

*(3) For* $\mathbb{P}$-*a.a.* $\omega$, *given* $\varepsilon > 0$, *there exists* $K = K(\omega, \varepsilon) < \infty$ *such that*

$$(2.8) \qquad |\chi(\omega, x)| \leq K + \varepsilon |x| \qquad \text{for all } x \in \xi.$$

*(4) For* $\mathbb{P}$-*a.a.* $\omega$, *the discrete-time process*

$$M(n) := \hat{X}^\omega(n) + \chi(\omega, \hat{X}^\omega(n))$$

*is a martingale (with respect to the natural* $\hat{X}^\omega$-*filtration) with square-integrable increments.*

We stress that the property $\mathbb{E}[(c_{0,x_1})^{-1}] < \infty$ is used to derive property (2).

PROOF OF LEMMA 2.1. (1) follows from the definition of the corrector $\chi$. To prove (2) note that there exists a sequence $g_n$ of bounded functions on $\mathcal{N}_0$ such that $\nabla g_n \to \chi$ in $\mathcal{H}$. From stationarity $\mathbb{E}[\nabla_{x_1} g_n] = 0$, so that the claim for $i = 1$ follows from the Schwarz inequality:

$$\mathbb{E}[|\chi(\omega, x_1) - \nabla_{x_1} g_n|] \leq \mathbb{E}[c_{0,x_1}^{-1}]^{1/2} \mathbb{E}[c_{0,x_1} |\chi(\omega, x_1) - \nabla_{x_1} g_n|^2]^{1/2}$$
$$\leq \mathbb{E}[c_{0,x_1}^{-1}]^{1/2} \|\chi - \nabla g_n\|_2 \to 0.$$

For any $i \in \mathbb{Z}$ we have, by the co-cycle property

$$\chi(\omega, x_i) = \chi(\omega, x_1) + \cdots + \chi(\tau_{x_{i-1}} \omega, x_i - x_{i-1}).$$

Therefore the claim follows for any $i$ by stationarity.

To prove (3), we shall use property (2). Namely, for $i = 1$ we have $\mathbb{E}[\chi(\omega, x_1)] = 0$. From the ergodic theorem and the co-cycle property we have

$$\frac{1}{x_n}\chi(\omega, x_n) = \frac{n}{x_n}\frac{1}{n}\sum_{j=0}^{n-1}[\chi(\omega, x_{j+1}) - \chi(\omega, x_j)]$$

$$= \frac{n}{x_n}\frac{1}{n}\sum_{j=0}^{n-1}[\chi(\tau_{x_j}\omega, x_{j+1} - x_j)] \to \mathbb{E}[Z_0]^{-1}\mathbb{E}[\chi(\omega, x_1)] = 0,$$

almost surely as $n \to \infty$. Similarly, one has $\frac{1}{x_{-n}}\chi(\omega, x_{-n}) \to 0$, $n \to \infty$. This proves the sublinearity of the corrector claimed in (2.8).

Finally, let us prove (4). Since $\hat{X}^\omega(n) = \sum_{j=0}^{n-1} f(\omega_j, \omega_{j+1})$, due to (2.2) and the stationarity of the process $\omega_n$ we know that $E_{\mathbb{Q}}[\hat{E}^\omega(\hat{X}^\omega(n)^2)] < \infty$, where $\hat{E}^\omega$ denotes the expectation with respect to the law of $\hat{X}^\omega$. Hence for $\mathbb{P}$-a.a. $\omega$ the random walk $\hat{X}^\omega(n)$ has square-integrable increments. Due to the sublinearity of the corrector $\chi$ given by (2.8), the same holds for the Markov chain $\chi(\omega, \hat{X}^\omega(n))$. Hence, for $\mathbb{P}$-a.a. $\omega$, $M(n)$ has square-integrable increments. The martingale property follows from the fact that for $\mathbb{P}$-a.a. $\omega$ the map $\xi \ni x \to x + \chi(\omega, x) \in \mathbb{R}$ is harmonic for the random walk $\hat{X}^\omega$ (see the discussion preceding the lemma). $\square$

We can now conclude the proof of the a.s. FCLT. Following line by line the $d = 2$ argument in [3], Section 6.1, one obtains an a.s. FCLT for the deformed walk $M(n) = \hat{X}^\omega(n) + \chi(\omega, \hat{X}^\omega(n))$. From this result and the sublinearity of the corrector field [see Lemma 2.1(3)], reasoning as in [3], Section 6.2, (6.10)–(6.13), one derives the FCLT for $\hat{X}^\omega$ for a suitable diffusion coefficient $\hat{D}(\beta)$. This coefficient must coincide with the one obtained in the proof of part (i). Finally, one derives the FCLT for $X^\omega$ by a time–change argument as in [3], Section 6.3, proof of Theorem 1.2.

**3. Proof of Theorem 1.2: Lower bounds on $D(\beta)$ by comparison with nearest-neighbor random walks.** In this section we prove the lower bounds (1.15) and (1.19) (left-hand side) by comparing Mott variable-range random walk with suitable nearest-neighbor random walks, whose diffusion coefficient can be computed explicitly. Trivially, the lower bound in (1.17) is implied by (1.15).

3.1. *Nearest-neighbor walks on $\xi$.* Let $\omega = (\xi, \Theta)$ denote a given realization of our marked point process [as defined in the Introduction; cf. the paragraph after (1.7)].

Given $\kappa \in (0, \infty]$, we consider the nearest-neighbor random walk $(Y_\kappa^\omega(t), t \geq 0)$ on $\xi$ starting at 0, with infinitesimal generator

$$\mathcal{L}_\kappa^\omega f(x_i) = \sum_{j : j = i \pm 1} r_{x_i, x_j}^{(\kappa)}(\omega)(f(x_j) - f(x_i)) \qquad \forall f : \xi \to \mathbb{R},$$

where

$$r_{x_i, x_j}^{(\kappa)}(\omega) = \exp\{-(|x_i - x_j|^\alpha \wedge \kappa) - \beta u(E_i, E_j)\}.$$

The random walk $Y_\kappa^\omega$ can be described as follows: after arriving at site $x_i$ the particle waits an exponential time with parameter

$$(3.1) \qquad \gamma_{x_i}^{(\kappa)}(\omega) = r_{x_i, x_{i+1}}^{(\kappa)}(\omega) + r_{x_i, x_{i-1}}^{(\kappa)}(\omega)$$

and then jumps to the nearest-neighbor point $x_{i \pm 1}$ with probability $r_{x_i, x_{i \pm 1}}^{(\kappa)}(\omega) / \gamma_{x_i}^{(\kappa)}(\omega)$. Since the parameters (3.1) satisfy $0 < \gamma_{x_i}^{(\kappa)}(\omega) \leq 2$, it is standard to check that the random walk is well defined (see, e.g., [5]) for any $0 < \kappa \leq \infty$. Since

$$(3.2) \qquad r_{x_i, x_j}^{(\infty)} = c_{x_i, x_j},$$

$Y_\infty^\omega$ is the nearest-neighbor version of the random walk defined by (1.8). Note that the definition above is such that when $\kappa < \infty$ one has uniformly elliptic jump rates.

By the same methods used in Section 2 and thanks to the finite second moment condition, it is simple to establish the following invariance principle. Let $P_{\kappa, \varepsilon}^\omega$ denote the law on $D([0, \infty), \mathbb{R})$ of the rescaled process $\varepsilon Y_\kappa^\omega(\varepsilon^{-2} t)$. Then, for any $0 < \kappa \leq \infty$, $P_{\kappa, \varepsilon}^\omega$ weakly converges in $\mathbb{P}$-probability to a one-dimensional Brownian motion $W$ with variance $\mathbb{E}(W(t)^2) = D_\kappa(\beta)t$, where

$$(3.3) \qquad D_\kappa(\beta) = \inf_{g \in L^\infty(\mathbb{P})} \mathbb{E}\left\{ \sum_{i = \pm 1} r_{0, x_i}^{(\kappa)}(x_i + \nabla_{x_i} g(\omega))^2 \right\}.$$

Note that, at this stage, the value of $D_\kappa(\beta)$ in (3.3) could well be zero.

As before we shall obtain an almost sure invariance principle under the following extra assumption:

$$(3.4) \qquad \mathbb{E}[1/r_{0, Z_0}^{(\kappa)}] < \infty.$$

Note that (3.4) is trivially satisfied for $\kappa < \infty$. Below we derive an explicit expression for $D_\kappa(\beta)$.

PROPOSITION 3.1. _Assume (3.4). Then $\mathbb{P}$-a.s. the law of $Y_\varepsilon^\omega = (\varepsilon Y_\kappa^\omega(\varepsilon^{-2} t) : t \geq 0)$ converges weakly as $\varepsilon \to 0$ to the law of the 1D Brownian motion $W$ with variance $\mathbb{E}(W^2(t)) = \widetilde{D}_\kappa(\beta)t$, where_

$$(3.5) \qquad \widetilde{D}_\kappa(\beta) = \frac{2\mathbb{E}[Z_0]^2}{\mathbb{E}(1/r_{0, x_1}^{(\kappa)})}.$$

The proof of the above proposition will be based on the corrector field. Alternatively, the same result could be obtained by expressing the nearest-neighbor random walk as a space–time change of a Brownian motion and applying Stone's method [12, 24].

PROOF. We repeat the main steps of the proof of Theorem 1.1 with the explicit choice for $\chi$ given by

$$
\begin{aligned}
(3.6) \qquad \chi(\omega, x_n) &= \sum_{j=0}^{n-1} \Big( \frac{C}{r_{x_j, x_{j+1}}^{(\kappa)}} - Z_j \Big), \\
\chi(\omega, x_{-n}) &= -\sum_{j=1}^{n} \Big( \frac{C}{r_{x_{-j}, x_{-j+1}}^{(\kappa)}} - Z_{-j} \Big),
\end{aligned}
$$

where

$$
C := \frac{\mathbb{E}[Z_0]}{\mathbb{E}(1/r_{0, x_1}^{(\kappa)})}.
$$

The crucial observation is that this function $\chi$ has the property (1), (2), (3) and (4) listed in Lemma 2.1. For instance, property (4) follows from the harmonicity $\mathcal{L}_\kappa^\omega \varphi(\omega, x) = 0$, $x \in \xi$, where $\varphi(\omega, x) = \chi(\omega, x) + x$. Therefore, it satisfies the same conclusions. In particular, the almost sure invariance principle holds for $\varepsilon Y_\kappa^\omega(\varepsilon^{-2} t)$ with constant diffusion coefficient given by $\widetilde{D}_\kappa(\beta) = \frac{1}{t} \mathbb{E} E_{0,\omega}[M_t^2]$, for any $t > 0$, where $M_t$ denotes the continuous-time martingale $M_t = \varphi(\omega, Y_\kappa^\omega(t))$. Here $E_{0,\omega}$ denotes expectation with respect to the random walk with generator $\mathcal{L}_\kappa^\omega$, started at the origin. It remains to check that $\widetilde{D}_\kappa(\beta)$ is given by (3.5). There are several ways to do this. For instance one can check that, $[\mathcal{L}_\kappa^\omega \varphi^2](\omega, 0) = C^2 / r_{0, x_1}^{(\kappa)} + C^2 / r_{x_{-1}, 0}^{(\kappa)}$. In this way it follows that

$$
(3.7) \quad \widetilde{D}_\kappa(\beta) = \lim_{t \to 0} \frac{1}{t} \mathbb{E} E_{0,\omega}[M_t^2] = \mathbb{E}[(\mathcal{L}_\kappa^\omega \varphi^2)(\omega, 0)] = 2C^2 \mathbb{E}[1/r_{0, x_1}^{(\kappa)}],
$$

which gives the desired expression (see, e.g., [8], (4.22), for a similar argument). $\square$

As a corollary we thus obtain the following consequences for the constant $D_\kappa(\beta)$ appearing in (3.3).

COROLLARY 3.2.   *For any $0 < \kappa \leq \infty$ we have*

$$
(3.8) \qquad D_\kappa(\beta) = \frac{2\mathbb{E}[Z_0]^2}{\mathbb{E}(1/r_{0, x_1}^{(\kappa)})}.
$$

*In particular, $D_\infty(\beta) > 0$ if and only if*

$$\mathbb{E}(\exp(Z_0^\alpha)) < \infty. \tag{3.9}$$

*Moreover, assuming (3.9), there exist $C_1, \kappa_1 \in (0, \infty)$ such that for all $\beta \geq 0$*

$$D_\infty(\beta) \geq C_1 e^{-\kappa_1 \beta}. \tag{3.10}$$

*Finally, if $\nu \neq \delta_0$ then there exist $C_2, \kappa_2 \in (0, \infty)$ such that for all $\beta \geq 0$*

$$D_\infty(\beta) \leq C_2 e^{-\kappa_2 \beta}. \tag{3.11}$$

PROOF. Clearly, whenever (3.4) is satisfied then $D_\kappa(\beta)$ must equal $\widetilde{D}_\kappa(\beta)$ and therefore (3.8) holds in this case. The only way to violate (3.4) is to have $\kappa = \infty$ and $\mathbb{E}(\exp(Z_0^\alpha)) = \infty$. But here we can use monotonicity in $\kappa$ of the variational characterization (3.3) which implies that

$$D_\infty(\beta) \leq D_\kappa(\beta) = \frac{2\mathbb{E}[Z_0]^2}{\mathbb{E}(1/r_{0,x_1}^{(\kappa)})} \qquad \forall \kappa \in (0, \infty).$$

By taking the limit $\kappa \to \infty$ we get that $D_\infty(\beta) = 0$. In particular, (3.8) is always valid.

To prove (3.10) simply observe that by (1.9), using $|E_i| \leq 1$ we have

$$\nu[\exp(\beta u(E_0, E_1))] \leq e^{2c_2\beta}. \tag{3.12}$$

Finally, if $\nu \neq \delta_0$ then there is $\varepsilon > 0$ such that $\nu(|E_0| \geq \varepsilon) > 0$, so that by (1.9)

$$\nu[\exp(\beta u(E_0, E_1))] \geq \nu(|E_0| \geq \varepsilon)e^{c_1\varepsilon\beta}, \tag{3.13}$$

which proves (3.11). □

3.2. *Proof of the lower bound (1.15).* We compare the variational characterizations (1.14) and (3.3). Clearly,

$$D(\beta) \geq D_\infty(\beta). \tag{3.14}$$

Hence the lower bound on $D(\beta)$ follows from (3.10).

3.3. *Proof of the lower bound (1.19).* Here $\{Z_j\}$ are i.i.d. with $\mathbb{E}[\exp \varepsilon Z_j] < \infty$ for some parameter $\varepsilon > 0$. Let $E_* > 0$ be a small fixed number [later we shall take $E_* = E_*(\beta) \to 0$ as $\beta \to \infty$]. From (1.14) we have

$$D(\beta) \geq \inf_{g \in L^\infty(\mathbb{P})} \mathbb{E}\left\{ \chi_{\{|E_0| \leq E_*\}} \sum_{i \in \mathbb{Z}} c_{0,x_i}(x_i + \nabla_{x_i}g)^2 \chi_{\{|E_i| \leq E_*\}} \right\}, \tag{3.15}$$

where $\chi$ denotes the indicator function. Let $\mathbb{P}_{0,*} = \mathbb{P}(\cdot \mid |E_0| \leq E_*)$ and let $\nu_* = \nu([-E_*, E_*])$, so that the right-hand side above can be written as

$$(3.16) \qquad \nu_* \inf_{g \in L^\infty(\mathbb{P})} \mathbb{E}_{0,*} \left\{ \sum_{i \in \mathbb{Z}} c_{0,x_i}(x_i + \nabla_i g)^2 \chi_{\{|E_i| \leq E_*\}} \right\},$$

where $\mathbb{E}_{0,*}$ denotes expectation with respect to $\mathbb{P}_{0,*}$. Let $\xi^* := \{x_i \in \xi : |E_i| \leq E_*\}$ denote the set of points in $\xi$ with energy of modulus less than $E_*$. This is a thinning of the original process. Moreover, the law $\mathbb{P}_*$ of $\xi^*$ under $\mathbb{P}_{0,*}$ coincides with the law of the process obtained from (1.7) where the distances $Z_j$ between consecutive points are replaced by new independent variables $Z_j^*$, each distributed as

$$(3.17) \qquad Z_1^* = \sum_{j=1}^Q Z_j,$$

where $Q$ is an independent geometric random variable with parameter $p := \nu_*$. That is, $Q$ is independent of the $\{Z_i\}$ and $\mathbb{P}(Q = k) = \nu_*(1 - \nu_*)^{k-1}$ for all positive integers $k$.

We write $\hat{c}_{x_i, x_j}$ for the rates (1.8) evaluated at $\beta = 0$, that is, $\hat{c}_{x_i, x_j} = \exp\{-|x_i - x_j|^\alpha\}$. In particular, by (1.9), assuming $|E_0| \leq E_*$ we have

$$c_{0,x_i} \chi_{\{|E_i| \leq E_*\}} \geq e^{-c\beta E_*} \hat{c}_{0,x_i} \chi_{\{|E_i| \leq E_*\}}$$

for some constant $c > 0$, for every $x_i$. From (3.15), (3.16) and by the same arguments used in the proof of Proposition 5 in [14], it then follows that

$$(3.18) \qquad D(\beta) \geq \nu_* e^{-c\beta E_*} D^*,$$

where $D^*$ is defined by

$$(3.19) \qquad D^* := \inf_{g \in L^\infty(\mathbb{P}_*)} \mathbb{E}_* \left\{ \sum_{i \in \mathbb{Z}} \hat{c}_{0,x_i}(x_i + \nabla_{x_i} g)^2 \right\}$$

with $\mathbb{E}_*$ denoting expectation with respect to the new measure $\mathbb{P}_*$. Thus $D^*$ is the diffusion coefficient at $\beta = 0$ associated to the point process $\mathbb{P}_*$. Now, we note that $D^*$ can be bounded from below by comparison with the associated nearest-neighbor walk as we did in (3.14). Namely,

$$(3.20) \qquad D^* \geq D_\infty^* := \frac{2\mathbb{E}[Z_1^*]^2}{\mathbb{E}(\exp[(Z_1^*)^\alpha])},$$

where $Z_1^*$ is defined by (3.17).

Next, we claim that there exists a constant $C$ such that for all values of $\nu_*$:

$$(3.21) \qquad \mathbb{E}[\exp[(Z_1^*)^\alpha]] \leq \exp[C\nu_*^{-\alpha/(1-\alpha)}].$$

To prove our claim we need the following lemma.

LEMMA 3.3. *Set* $\varphi_*(t) := \mathbb{P}(Z_1^* \geq t)$. *Then there exists a constant* $c > 0$ *independent of* $\nu_*$ *such that for all* $t \geq 0$

$$\varphi_*(t) \leq 4e^{-c\nu_* t}. \tag{3.22}$$

PROOF. Since $\varphi_*$ is bounded, it is enough to prove the statement for $\nu_*$ small enough. Let $\Lambda(\delta) := \log \mathbb{E}[e^{\delta Z_1}]$, $\delta \geq 0$, denote the logarithmic moment generating function of the original variable $Z_1$. We know this is finite for sufficiently small $\delta > 0$ since we are assuming $\mathbb{E}[e^{\varepsilon Z_1}] < \infty$ for some $\varepsilon > 0$. From the independence of the $\{Z_i\}$ we have $\mathbb{E}[e^{\delta S_n}] = e^{\Lambda(\delta)n}$, where $S_n = \sum_{i=1}^n Z_i$. Since the geometric variable $Q$ is independent of the $\{Z_i\}$ we have

$$\mathbb{E}[e^{\delta Z_1^*}] = \mathbb{E}[e^{\delta S_Q}] = \mathbb{E}[e^{\Lambda(\delta)Q}].$$

Moreover, for any $a > 0$ such that $e^a < (1 - \nu_*)^{-1}$, for the geometric variable $Q$ with parameter $\nu_*$ we have

$$\mathbb{E}[e^{aQ}] = \frac{\nu_* e^a}{1 - (1 - \nu_*)e^a}.$$

Therefore, for any $\delta > 0$ such that $e^{\Lambda(\delta)} < (1 - \nu_*)^{-1}$ we have

$$\varphi_*(t) \leq e^{-\delta t} \mathbb{E}[e^{\delta Z_1^*}] \leq e^{-\delta t} \frac{\nu_* e^{\Lambda(\delta)}}{1 - (1 - \nu_*)e^{\Lambda(\delta)}}.$$

It follows that for any $\delta > 0$ such that

$$e^{\Lambda(\delta)} \leq \frac{1 - \nu_*/2}{1 - \nu_*}, \tag{3.23}$$

we have

$$\varphi_*(t) \leq 2e^{\Lambda(\delta)}e^{-\delta t}. \tag{3.24}$$

Let us now show that there exists a constant $c > 0$ such that setting $\delta = c\nu_*$ we can satisfy the bound (3.23) for all sufficiently small values of $\nu_*$. In this case the desired estimate (3.22) would follow from (3.24) since $\Lambda(\delta) \to 1$ as $\delta \to 0$.

Set $\mu := \mathbb{E}[Z_1] > 0$, so that $\Lambda'(0) = \mu$. It is easy to check that $\Lambda(\delta) \leq 2\mu\delta$ for all sufficiently small values of $\delta$. Therefore (3.23) follows if

$$\delta \leq \frac{1}{2\mu} \log \frac{(1 - \nu_*/2)}{1 - \nu_*}.$$

However, one checks that $\log \frac{(1 - \nu_*/2)}{1 - \nu_*} \geq \nu_*/2$, and therefore the estimate above is certainly satisfied by $\delta = \frac{1}{4\mu}\nu_*$. $\square$

3.3.1. *Proof of claim (3.21).* Using Lemma 3.3, we have

$$
\begin{aligned}
\mathbb{E}(\exp[(Z_1^*)^\alpha]) &= \int_0^\infty \mathbb{P}(\exp[(Z_1^*)^\alpha] \geq t)\, dt \\
&= 1 + \int_1^\infty \varphi_*[(\log t)^{1/\alpha}]\, dt \\
&\leq 1 + 4 \int_1^\infty \exp[-c\nu_*(\log t)^{1/\alpha}]\, dt.
\end{aligned}
\tag{3.25}
$$

Since $\alpha < 1$, the last integral is finite. In particular, it is enough to prove (3.21) for $\nu_*$ small enough. By the variable change $s = \lambda_*(\log t)^{1/\alpha}$, where $\lambda_* := c\nu_*$, the last integral becomes

$$
\lambda_*^{-\alpha} \int_0^\infty \alpha s^{\alpha-1} e^{-s} e^{s^\alpha \lambda_*^{-\alpha}}\, ds.
$$

This, in turn, is estimated from above by

$$
C_\alpha \lambda_*^{-\alpha} e^{\lambda_*^{-\alpha}} + \lambda_*^{-\alpha} \int_1^\infty e^{-s} e^{s^\alpha \lambda_*^{-\alpha}}\, ds,
$$

with some constant $C_\alpha$ independent of $\lambda_*$. To estimate the last integral we set $C = 2^{1/(1-\alpha)}$ and observe that for $s > C\lambda_*^{-\alpha/(1-\alpha)}$ one has $s^\alpha \lambda_*^{-\alpha} < s/2$. Therefore, the integrand $e^{-s} e^{s^\alpha \lambda_*^{-\alpha}}$ can be estimated by $\exp[C^\alpha \lambda_*^{-\alpha/(1-\alpha)}]$ for $s \leq C\lambda_*^{-\alpha/(1-\alpha)}$ and by $e^{-s/2}$ for $s > C\lambda_*^{-\alpha/(1-\alpha)}$. By splitting the integral one obtains

$$
\int_1^\infty e^{-s} e^{s^\alpha \lambda_*^{-\alpha}}\, ds \leq C\lambda_*^{-\alpha/(1-\alpha)} \exp[C^\alpha \lambda_*^{-\alpha/(1-\alpha)}] + 2 \leq \exp[2C^\alpha \lambda_*^{-\alpha/(1-\alpha)}],
$$

where the last bound follows if we assume that $\lambda_*$ is sufficiently small. Putting these estimates together and using $\lambda_*^{-\alpha} e^{\lambda_*^{-\alpha}} \leq \exp[\lambda_*^{-\alpha/(1-\alpha)}]$ (for $\lambda_*$ small enough) we have

$$
C_\alpha \lambda_*^{-\alpha} e^{\lambda_*^{-\alpha}} + \lambda_*^{-\alpha} \int_1^\infty e^{-s} e^{s^\alpha \lambda_*^{-\alpha}}\, ds \leq \exp[C' \lambda_*^{-\alpha/(1-\alpha)}]
\tag{3.26}
$$

for some new constant $C'$, whenever $\lambda_*$ is sufficiently small. Using (3.25), and recalling that $\lambda_* = c\nu_*$, this proves the claim (3.21).

We can now finish the proof of the lower bound (1.19) in Theorem 1.2. Going back to (3.20) and observing that $\mathbb{E}[Z_1^*] \geq \mathbb{E}[Z_1] =: \mu > 0$ for some constant $\mu$, we see that thanks to (3.21) we have

$$
D_\infty^* \geq 2\mu^2 \exp[-C\nu_*^{-\alpha/(1-\alpha)}].
$$

From (3.18) and (3.20) it follows that

$$
D(\beta) \geq 2\mu^2 \nu^* \exp[-c\beta E_* - C\nu_*^{-\alpha/(1-\alpha)}].
\tag{3.27}
$$

By assumption there are positive constants $c_1, c_2$ such that

$$c_1 E_*^\delta \le \nu_* \le c_2 E_*^\delta. \tag{3.28}$$

Therefore, the exponent in (3.27) is bounded from below by $-c\beta E_* - \kappa_1 \times E_*^{-\delta\alpha/(1-\alpha)}$, with a new constant $\kappa_1$. Choosing $E_* = \beta^{-(1-\alpha)/(1-\alpha+\delta\alpha)}$ one obtains

$$D(\beta) \ge \exp[-\kappa \beta^{\delta\alpha/(1-\alpha+\delta\alpha)}]$$

for all sufficiently large $\beta$, for some constant $\kappa > 0$. This ends the proof of the lower bound (1.19) in Theorem 1.2.

We remark that the preceding proof is obtained by optimization over several possible choices of nearest-neighbor walks. In particular, it shows that the true diffusion constant $D(\beta)$ can be much larger than the diffusion constant of the nearest-neighbor walk defined in (3.2), which is characterized by the low-temperature behavior given by (3.10) and (3.11).

## 4. Proof of Theorem 1.2: subdiffusive regime.

Here we shall prove points (2) and (3)(i) of Theorem 1.2. Due to point (1) in Theorem 1.2, we only need to show for renewal point processes $\xi$ that $D(\beta) = 0$ if (i) $\alpha \ge 1$ and $\mathbb{E}[\exp(Z_j^\alpha)] = \infty$, or if (ii) $\alpha < 1$ and $\mathbb{E}[\exp(\gamma Z_j^\alpha)] = \infty$ for some $\gamma \in (0,1)$. To this end, we define $\varphi(x) = e^{|x|^\alpha}$ and, for any $g \in L^\infty(\mathbb{P})$, we set

$$A(g) = A_+(g) + A_-(g), \qquad A_\pm(g) := \mathbb{E}\left\{ \sum_{i=1}^\infty \varphi(x_{\pm i})(x_{\pm i} + \nabla_{x_{\pm i}} g)^2 \right\}. \tag{4.1}$$

Then, due to (1.9), we have the upper bound $D(\beta) \le c(\beta) \inf_{g \in L^\infty(\mathbb{P})} A(g)$. We only need to show that this last infimum is zero in the above cases (i) and (ii).

To this end we pick $g \in L^\infty(\mathbb{P})$ and, given $i \ge 1$, we write

$$x_i + \nabla_{x_i} g = \sum_{\ell=0}^{i-1} [Z_\ell + \nabla_{Z_\ell} g(\tau_{x_\ell} \omega)].$$

Applying the Schwarz inequality we get

$$(x_i + \nabla_{x_i} g)^2 \le i \sum_{\ell=0}^{i-1} [Z_\ell + \nabla_{Z_\ell} g(\tau_{x_\ell} \omega)]^2.$$

Therefore, we can bound

$$\begin{aligned}
A_+(g) &\le \mathbb{E}\left\{ \sum_{i=1}^\infty i\varphi(x_i) \sum_{\ell=0}^{i-1} [Z_\ell + \nabla_{Z_\ell} g(\tau_{x_\ell} \omega)]^2 \right\} \\
&= \sum_{\ell=0}^\infty \mathbb{E}\{c_\ell(\xi)[Z_\ell + \nabla_{Z_\ell} g(\tau_{x_\ell} \omega)]^2\},
\end{aligned} \tag{4.2}$$

where we use the notation $c_\ell(\xi) := \sum_{k=\ell+1}^\infty k\varphi(x_k)$. Using shift invariance, we can write

$$\mathbb{E}\{c_\ell(\xi)[Z_\ell + \nabla_{Z_\ell} g(\tau_{x_\ell}\omega)]^2\} = \mathbb{E}\left\{\sum_{k=1}^\infty ((\ell+k)\varphi(x_k - x_{-\ell}))[Z_0 + \nabla_{Z_0} g(\omega)]^2\right\}.$$

In conclusion, we arrive at the bound $A_+(g) \leq \mathbb{E}\{C_1(\xi)[Z_0 + \nabla_{Z_0} g(\omega)]^2\}$, where

$$C_1(\xi) := \sum_{m=0}^\infty \sum_{n=0}^\infty (1+n+m)\varphi\left(Z_0 + \sum_{\ell=1}^n Z_\ell + \sum_{\ell=1}^m Z_{-\ell}\right),$$

(we agree to set $\sum_{\ell=a}^b$ equal to zero if $a > b$). By the same arguments, one can bound $A_-(g)$ and therefore arrive at the estimate

$$A(g) \leq \mathbb{E}\left\{\sum_{j=\pm 1} C_j(\xi)[x_j + \nabla_{x_j} g(\omega)]^2\right\},$$

where

$$C_{-1}(\xi) := \sum_{m=0}^\infty \sum_{n=0}^\infty (1+n+m)\varphi\left(Z_{-1} + \sum_{\ell=0}^{n-1} Z_\ell + \sum_{\ell=2}^{m+1} Z_{-\ell}\right).$$

Collecting all our bounds, we get that

$$(4.3) \qquad D(\beta) \leq c(\beta) \inf_{g \in L^\infty(\mathbb{P})} A(g) \leq c(\beta)\bar{D},$$

where

$$(4.4) \qquad \bar{D} := \inf_{\substack{g=g(\xi) \\ \text{bounded}}} \mathbb{E}\left\{\sum_{j=\pm 1} C_j(\xi)[x_j + \nabla_{x_j} g(\xi)]^2\right\}.$$

Given $\kappa > 0$ we define

$$C_j^{(\kappa)}(\xi) := C_j(\xi) + \kappa.$$

Since $C_j(\xi) \leq C_j^{(\kappa)}(\xi)$, we conclude that

$$\bar{D} \leq \bar{D}_\kappa := \inf_{\substack{g=g(\xi) \\ \text{bounded}}} \mathbb{E}\left\{\sum_{j=\pm 1} C_j^{(\kappa)}(\xi)[x_j + \nabla_{x_j} g(\xi)]^2\right\}.$$

Since $\mathbb{E}(\sum_{j=\pm 1} C_j^{(\kappa)}(\xi)x_j^2) < \infty$ [recall that $\mathbb{E}(Z_j^2) < \infty$] and $\mathbb{E}(1/C_j^{(\kappa)}) < \infty$ we can apply the same arguments used in Section 3 and check that $\bar{D}_\kappa$ is the diffusion coefficient of the nearest-neighbor random walk on $\xi$ jumping from $x_i$ to $x_{i\pm 1}$ with rate $C_{\pm 1}^{(\kappa)}(\tau_{x_i}\xi)$, and that $\bar{D}_\kappa$ equals $2\mathbb{E}(Z_0)^2/\mathbb{E}(1/C_1^{(\kappa)})$.

Since $C_1^{(\kappa)}(\xi) \searrow C_1(\xi)$, by applying the monotone convergence theorem we conclude that

$$\bar{D} \leq \lim_{\kappa \downarrow \infty} 2\mathbb{E}(Z_0)^2/\mathbb{E}(1/C_1^{(\kappa)}) = 2\mathbb{E}(Z_0)^2/\mathbb{E}(1/C_1).$$

In conclusion, in order to prove that $D(\beta) = 0$ it is enough to show that $\mathbb{E}(1/C_1) = \infty$.

4.1. *Case $\alpha \geq 1$ and $\mathbb{E}[\exp(Z_j^\alpha)] = \infty$.* Since $\alpha \geq 1$, we have

$$\varphi\left(Z_0 + \sum_{\ell=1}^n Z_\ell + \sum_{\ell=1}^m Z_{-\ell}\right) \leq \varphi(Z_0)\varphi\left(\sum_{\ell=1}^n Z_\ell + \sum_{\ell=1}^m Z_{-\ell}\right).$$

It follows that $C_1(\xi) \leq \varphi(Z_0)\bar{C}_1(\xi)$, where

$$\bar{C}_1(\xi) := \sum_{m=0}^\infty \sum_{n=0}^\infty (1 + n + m)\varphi\left(\sum_{\ell=1}^n Z_\ell + \sum_{\ell=1}^m Z_{-\ell}\right).$$

It is easy to check that $\mathbb{E}[(\bar{C}_1)^{-1}]^{-1} \leq \mathbb{E}[\bar{C}_1] \leq K$ for some constant $K < \infty$. Moreover, by independence of the $Z_k$

$$\mathbb{E}[(C_1(\xi))^{-1}] \geq \mathbb{E}[\varphi(Z_0)^{-1}]\mathbb{E}[(\bar{C}_1(\xi))^{-1}] \geq K^{-1}\mathbb{E}[\varphi(Z_0)^{-1}] = \infty,$$

due to our main assumption.

4.2. *Case $\alpha < 1$ and $\mathbb{E}[\exp(\gamma Z_j^\alpha)] = \infty$ for some $\gamma \in (0,1)$.* It is the same proof as above, with the exception that by concavity of $x^\alpha$ we can estimate

$$(a + b)^\alpha \geq (\gamma a + (1 - \gamma)b)^\alpha \geq \gamma a^\alpha + (1 - \gamma)b^\alpha$$

for all $\gamma \in (0,1)$ and all positive $a, b$, therefore

$$\varphi\left(Z_0 + \sum_{\ell=1}^n Z_\ell + \sum_{\ell=1}^m Z_{-\ell}\right) \leq \exp(-\gamma Z_0^\alpha)\exp\left((1-\gamma)\left[\sum_{\ell=1}^n Z_\ell + \sum_{\ell=1}^m Z_{-\ell}\right]^\alpha\right).$$

From here on we can repeat the previous argument and the conclusion follows.

## 5. Proof of Theorem 1.2: Upper bounds on $D(\beta)$.

In this section we prove the upper bounds in (1.17) and (1.19) by using suitable test functions $g$ in the variational formula (1.14). To this end the exponential moment assumption is not needed. Indeed, what we shall use is only the existence of the second moment in Lemma 5.1 below.

We start the proof in the general case $\alpha > 0$ and then separate the two cases $\alpha \geq 1$ and $\alpha < 1$ afterward. Given a value $E_* > 0$ of the energy and a realization $\omega = (\xi, \Theta)$ of the sequence $\{x_j\}$ in (1.7) and of the associated energies $E_j$, we define the new sequence $y_k$ as follows:

$$(5.1) \quad y_0 := 0, \qquad y_k := \begin{cases} \inf\{x_j : x_j > y_{k-1}, |E_j| \leq E_*\}, & k \geq 1, \\ \sup\{x_j : x_j < y_{k+1}, |E_j| \leq E_*\}, & k \leq -1. \end{cases}$$

Note that, apart from $y_0 = x_0 = 0$ which may have energy $E_0$ satisfying $|E_0| > E_*$, the sequence $y_k$ corresponds exactly to the points $x_j$ such that $|E_j| \leq E_*$. We call $\xi^* := \{y_k, k \in \mathbb{Z}\}$ the new point process (a thinning of $\xi$). Next, given a sequence of nonnegative numbers $\{c(n), n \in \mathbb{N}\}$ to be chosen later we define

$$(5.2) \qquad \begin{aligned} L_n^+ &:= \inf\{k \geq 0 : y_{k+1} - y_k \geq c(n)\}, \\ L_n^- &:= \sup\{k \leq 0 : y_k - y_{k-1} \geq c(n)\}. \end{aligned}$$

We consider the function $g(\omega)$ given by

$$(5.3) \qquad g(\omega) := y_{L_n^+} \wedge n.$$

Therefore, from (1.14) we may estimate

$$(5.4) \qquad D(\beta) \leq \mathbb{E}\left\{\sum_{j \in \mathbb{Z}} e^{-|x_j|^\alpha} e^{-\beta u(E_0, E_j)} (x_j + \nabla_{x_j} g)^2\right\}.$$

We split the expectation above in two terms given by the events $|E_0| \leq E_*$ and its complement. Letting $\mathbb{E}_{0,*}$ denote the expectation with respect to the conditional probability $\mathbb{P}_{0,*} = \mathbb{P}(\cdot | |E_0| \leq E_*)$ and setting $\nu_* = \nu([-E_*, E_*])$ we have

$$(5.5) \qquad D(\beta) = D_1(\beta) + D_2(\beta),$$

where

$$(5.6) \quad D_1(\beta) := \nu_* \mathbb{E}_{0,*}\left\{\sum_{j \in \mathbb{Z}} e^{-|x_j|^\alpha} e^{-\beta u(E_0, E_j)} (x_j + \nabla_{x_j} g)^2\right\},$$

$$(5.7) \quad D_2(\beta) := \mathbb{E}\left\{\sum_{j \in \mathbb{Z}} e^{-|x_j|^\alpha} e^{-\beta u(E_0, E_j)} (x_j + \nabla_{x_j} g)^2 \chi(|E_0| > E_*)\right\}.$$

We first estimate the term $D_2(\beta)$. Here we know that $|E_0| > E_*$ so that from (1.9) and the fact that for $\gamma < 1$ we can bound $e^{-|x|^\alpha}(x + \nabla_x g)^2 \leq cn^2 e^{-\gamma |x|^\alpha}$ for some finite $c = c(\gamma)$ we have

$$(5.8) \qquad D_2(\beta) \leq cn^2 e^{-c_1 \beta E_*} \mathbb{E}\left\{\sum_{j \in \mathbb{Z}} e^{-\gamma |x_j|^\alpha}\right\}.$$

Using Lemma A.1 we see that the last expectation in (5.8) is finite and therefore we obtain the following estimate for $D_2(\beta)$, for some finite constant $C$:

$$(5.9) \qquad D_2(\beta) \leq Cn^2 e^{-c_1 \beta E_*}.$$

We turn to the estimate of $D_1(\beta)$. In (5.6) we can further split the sum depending on whether $|E_j| > E_*$ or not, that is, whether $x_j \in \xi^*$ or not. Estimating as in (5.8) and (5.9) we have

$$(5.10) \quad \mathbb{E}_{0,*}\left\{\sum_{j \in \mathbb{Z}} e^{-|x_j|^\alpha} e^{-\beta u(E_0, E_j)}(x_j + \nabla_{x_j} g)^2 \chi(x_j \notin \xi^*)\right\} \leq C n^2 e^{-c_1 \beta E_*}.$$

Therefore

$$(5.11) \quad D_1(\beta) \leq C n^2 e^{-c_1 \beta E_*} + D_1^*,$$

where, recalling that $\xi^* = \{y_j\}$ and neglecting the energy contribution, we define

$$(5.12) \quad D_1^* := \mathbb{E}_{0,*}\left\{\sum_{j \in \mathbb{Z}} e^{-|y_j|^\alpha}(y_i + \nabla_{y_i} g)^2\right\}.$$

We observe that the law of $\xi^*$ under $\mathbb{P}_{0,*}$ is given by the renewal process such that the distance between consecutive points $Z_j^* := y_{j+1} - y_j$ is distributed as $\sum_{\ell=1}^Q Z_\ell$, where $Q$ is an independent geometric random variable with parameter $p = \nu_*$. That is, $Q$ is independent of the $\{Z_i\}$ and $\mathbb{P}(Q = k) = \nu_*(1 - \nu_*)^{k-1}$, $k \in \mathbb{N}$. Again for any given $0 < \gamma < 1$, we have $e^{-|x|^\alpha}(x + \nabla_x g)^2 \leq cn^2 e^{-\gamma|x|^\alpha}$ for some positive constant $c = c(\gamma)$. Therefore we can estimate

$$(5.13) \quad D_1^* \leq c(A_0 + A_1^+ + A_1^- + A_2^+ + A_2^-),$$

where

$$A_0 = \mathbb{E}_{0,*}\left\{\sum_{j \in \mathbb{Z}} e^{-|y_j|^\alpha}(y_j + \nabla_{y_j} g)^2 \chi(-n/2 \leq y_{L_n^-} \leq y_j \leq y_{L_n^+} \leq n/2)\right\},$$

$$A_1^+ = n^2 \mathbb{E}_{0,*}\left\{\sum_{j \in \mathbb{Z}} e^{-\gamma|y_j|^\alpha} \chi(y_{L_n^+} > n/2)\right\},$$

$$A_1^- = n^2 \mathbb{E}_{0,*}\left\{\sum_{j \in \mathbb{Z}} e^{-\gamma|y_j|^\alpha} \chi(y_{L_n^-} < -n/2)\right\},$$

$$A_2^+ = n^2 \mathbb{E}_{0,*}\left\{\sum_{j \in \mathbb{Z}} e^{-\gamma|y_j|^\alpha} \chi(y_j > y_{L_n^+})\right\},$$

$$A_2^- = n^2 \mathbb{E}_{0,*}\left\{\sum_{j \in \mathbb{Z}} e^{-\gamma|y_j|^\alpha} \chi(y_j < y_{L_n^-})\right\}.$$

The first observation is that $A_0 = 0$. Indeed, since $-n/2 \leq y_{L_n^-} \leq y_j \leq y_{L_n^+} \leq n/2$, we have that $g(\tau_{y_j} \xi^*) = y_{L_n^+} - y_j$ so that $y_j + \nabla_{y_j} g = 0$ in this case.

We turn to estimate $A_1^+$. Fix $p, q > 1$ with $1/p + 1/q = 1$ and set $\psi := \sum_{j \in \mathbb{Z}} e^{-\gamma |y_j|^\alpha}$. Then by Hölder's inequality

$$A_1^+ \leq n^2 \mathbb{P}_{0,*}(y_{L_n^+} > n/2)^{1/p} \mathbb{E}_{0,*}[\psi^q]^{1/q}.$$

Thanks to Lemma A.1 we have $\mathbb{E}_{0,*}[\psi^q]^{1/q} \leq C_q < \infty$ for any $q > 1$. Now, note that by construction $y_{L_n^+} \leq L_n^+ c(n)$ so that

$$\mathbb{P}_{0,*}(y_{L_n^+} > n/2) \leq \mathbb{P}_{0,*}(L_n^+ > n/(2c(n))).$$

Set $m_n = [n/(2c(n))]$ (the integer part). Then $L_n^+ > n/(2c(n))$ implies that the first $m_n$ variables $Z_j^*$, $j = 0, \ldots, m_n - 1$ are all smaller than $c(n)$, so that

$$(5.14) \quad \mathbb{P}_{0,*}(y_{L_n^+} > n/2) \leq (1 - \mathbb{P}(Z_1^* \geq c(n)))^{m_n} \leq \exp[-m_n \varphi_*(c(n))],$$

where $\varphi_*(t) := \mathbb{P}(Z_j^* \geq t)$. We have obtained

$$(5.15) \qquad\qquad A_1^+ \leq n^2 C_q \exp\left[-\frac{1}{p} m_n \varphi_*(c(n))\right].$$

Since $A_1^+ = A_1^-$ by symmetry, the same bound applies to $A_1^-$.

We now estimate the term $A_2^+ = A_2^-$. If $y_j \in \xi^*$ with $y_j > y_{L_n^+}$ then $y_j \geq y_{L_n^+} + c(n) + (y_j - y_{L_n^+ + 1}) =: a_1 + a_2 + a_3$. Since $a_i \geq 0$, we can use the elementary inequality $(\sum_{i=1}^3 a_i)^\alpha \geq d(\alpha) \sum_{i=1}^3 a_i^\alpha$, where $d(\alpha) := 1 \wedge 3^{\alpha-1}$. Hence, by the renewal property,

$$(5.16) \quad A_2^+ \leq n^2 e^{-\gamma d(\alpha) c(n)^\alpha} \mathbb{E}_{0,*}(\exp[-\gamma d(\alpha) y_{L_n^+}^\alpha]) \mathbb{E}_{0,*}\left(\sum_{j=0}^\infty e^{-\gamma d(\alpha) y_j^\alpha}\right).$$

Again Lemma A.1 allows to bound the last factor in the right-hand side by a finite constant $C$. Therefore, using the obvious bound $\mathbb{E}_{0,*}(\exp[-\gamma d(\alpha) y_{L_n^+}^\alpha]) \leq 1$ we have

$$(5.17) \qquad\qquad A_2^+ \leq C n^2 e^{-\gamma d(\alpha) c(n)^\alpha}.$$

Due to the above bounds and the arbitrariness of $p$ in (5.15), we obtain that, for all $\varepsilon > 0$ and $\gamma \in (0, 1)$, for a suitable positive constant $C$

$$(5.18) \quad D_1^* \leq C n^2 \{\exp[-(1-\varepsilon) m_n \varphi_*(c(n))] + \exp[-\gamma d(\alpha) c(n)^\alpha]\},$$

where $m_n = [n/(2c(n))]$, $d(\alpha) = 1 \wedge 3^{\alpha-1}$, and $\varphi_*(t) = \mathbb{P}(Z_1^* \geq t)$ denotes the tail probability of $Z_1^*$.

We need the following simple lemma.

LEMMA 5.1. *Suppose that $Z_1$ has finite second moment. There exists a constant $c > 0$ independent of $\nu_*$ such that*

$$\varphi_*(t) \geq c e^{-\lambda_* t}, \qquad t \geq 0,$$

*where $\lambda_* := -\frac{1}{\mu} \log(1 - \nu_*)$, $\mu := \mathbb{E}[Z_1]$.*

PROOF.   Set $S_n = \sum_{i=1}^{n} Z_i$. Then by the central limit theorem we have $\mathbb{P}(S_n \geq \mu n) \to \frac{1}{2}$, $n \to \infty$. Recall that $Q$ is an independent geometric random variable with parameter $p = \nu_*$. Therefore, for $t$ sufficiently large

$$\varphi(t) = \mathbb{P}(S_Q \geq t) \geq \mathbb{P}(S_Q \geq t; Q \geq t/\mu)$$

$$\geq \mathbb{P}(S_{[t/\mu]} \geq t)\mathbb{P}(Q \geq t/\mu) \geq c\mathbb{P}(Q \geq t/\mu) \geq ce^{-\lambda_* t}$$

with $\lambda(\nu_*) = -\frac{1}{\mu}\log(1 - \nu_*)$.   $\square$

We are now able to finish the proof of the upper bounds in (1.17) and (1.19). We start with the case $\alpha \geq 1$.

5.1. *Proof of the upper bound in (1.17).*   We choose the sequence $c(n)$ as

$$c(n) := a \log n$$

with $a > 0$. We want (5.18) to be smaller than an inverse power of $n$. Here $d(\alpha) = 1$. By Lemma 5.1 the first exponential in (5.18) is smaller than any inverse power of $n$ provided that $a\lambda_* < 1$. For the second exponential in (5.18) the worst case is when $\alpha = 1$. Here we need, for example, $a\gamma \geq 3$ to ensure that

$$n^2 \exp[-\gamma c(n)^\alpha] \leq \frac{1}{n}$$

for all $\alpha \geq 1$. Since $\gamma < 1$ is arbitrarily close to 1 we can reach this by requiring $a > 3$. Therefore we obtain, as soon as $a\lambda_* < 1$ and $a > 3$,

$$(5.19) \qquad D_1^* \leq \frac{2}{n}$$

for all sufficiently large $n$. Since we are assuming that $\nu_* \to 0$ as $E_* \to 0$ there is no difficulty in taking $E_*$ sufficiently small so that both $a > 3$ and $a\lambda_* < 1$ are satisfied. Finally, we can set

$$(5.20) \qquad n := e^{c_0 E_* \beta}$$

for a sufficiently small $c_0 > 0$, so that using (5.5), (5.9), (5.11) and (5.19) we obtain that $D(\beta) \leq e^{-\kappa\beta}$ for some $\kappa > 0$ and for all sufficiently large $\beta$. This ends the proof of the upper bound (1.17) in Theorem 1.2.

5.2. *Proof of the upper bound in (1.19).*   We are going to use the same idea we used for (1.17), with the difference that now the value $E_*$ will itself depend on $\beta$ and will go to 0 as $\beta \to \infty$. Recall that here we assume that $\nu_*$ is of order $E_*^\delta$ as $E_* \to 0$; cf. (3.28). We take $c(n)$ as

$$(5.21) \qquad c(n) := (a \log n)^{1/\alpha}.$$

As in the previous case we want an estimate as (5.19). Therefore the second exponential term in (5.18) requires the condition $\gamma 3^{\alpha-1} a \geq 3$. Since $\gamma < 1$ is arbitrarily close to 1 we can meet this requirement by taking $a > 3^{2-\alpha}$. For the first exponential, using Lemma 5.1, we see that

$$\lambda_* < (a \log n)^{-(1-\alpha)/\alpha},$$

is sufficient to have that $\exp[-(1-\varepsilon) m_n \varphi_*(c(n))]$ vanishes faster than any inverse power of $n$. Since $\lambda_* = -\mu^{-1} \log(1 - \nu_*) \leq 2\mu^{-1} \nu_*$ for small $\nu_*$, it is sufficient to have, for a suitable constant $C$,

$$(5.22) \qquad \nu_* < C (\log n)^{-(1-\alpha)/\alpha}.$$

The above argument shows that if we choose $c(n)$ as in (5.21) with, for example, $a = 9$ and suppose that $\nu_*$ satisfies (5.22) then we know that (5.19) holds for all $n$ sufficiently large.

Let us now rephrase things in terms of $E_*$ and $\beta$. We have $\nu_* \sim E_*^{\delta}$ and we may choose $n$ again as in (5.20). Therefore (5.22) can be rewritten, after suitably modifying the constant $C$,

$$(5.23) \qquad E_*^{\delta} < C (E_* \beta)^{-(1-\alpha)/\alpha}.$$

We can then choose $E_*(\beta) := c \beta^{-(1-\alpha)/(\delta\alpha+1-\alpha)}$ for a sufficiently small constant $c > 0$ and we are sure to satisfy (5.23). Once we have (5.23) we know that (5.19) holds with $n = e^{c_0 \beta E_*(\beta)}$, provided $\beta$ is large enough. Collecting the previous estimates in (5.5), (5.9) and (5.11) we have thus obtained, for suitable $\kappa_1, \kappa_2 > 0$:

$$D(\beta) \leq e^{-\kappa_1 \beta E_*(\beta)} = e^{-\kappa_2 \beta^{\delta\alpha/(\delta\alpha+1-\alpha)}}.$$

This ends the proof of Theorem 1.2.

**6. Bounds on Cheeger's constant and spectral gap.** In this section we are going to prove Proposition 1.3, Theorem 1.4 and Corollary 1.5. We recall that, as a convention, whenever we state that some event $E_L$ involving the parameter $L$ holds $\mathbb{P}$-a.s. we mean that with $\mathbb{P}$-probability 1 there exists $L_0 = L_0(\xi) < \infty$ such that the event $E_L$ occurs for all $L \in \mathbb{N}$, with $L \geq L_0$.

Let $\xi = \{x_k, k \in \mathbb{Z}\}$ denote as usual the renewal process generated by (1.7) and recall that $\xi_L := \xi \cap \Lambda_L$, $\Lambda_L = [-L/2, L/2]$. Recall the definition of Cheeger's isoperimetric constant (1.21) and the definition of $\zeta_L$ (1.22).

6.1. *Proof of Proposition 1.3.* To prove the left bound in (1.23), observe that for any nonempty subset $U$ with $\#(U) \leq \frac{1}{2} \#(\xi_L)$ we have

$$I_U(\xi) \geq \frac{2}{\#(\xi_L)} e^{-\zeta_L^{\alpha}}.$$

By ergodicity we have $\#(\xi_L)/L \to 1/\mathbb{E}(Z_i)$, $\mathbb{P}$-a.s. Therefore we may assume that

$$\frac{1}{C_1}L \leq \#(\xi_L) \leq C_1 L \tag{6.1}$$

for some constant $C_1$, and the claimed lower bound follows.

We turn to the proof of the upper bound in (1.23). We fix $\gamma \in (0,1)$ and let $a_{\gamma L} < b_{\gamma L}$ be two consecutive points realizing the maximum in (1.22) when $L$ is replaced by $\gamma L$, that is, $b_{\gamma L} - a_{\gamma L} = \zeta_{\gamma L}$. Since $\gamma < 1$, in analogy with (6.1) we know that taking $U$ to be either $U_1 := \xi \cap [-\frac{L}{2}, a_{\gamma L}]$ or $U_2 := \xi \cap [b_{\gamma L}, \frac{L}{2}]$ we have

$$c_\gamma L \leq \#(U) \leq \tfrac{1}{2}\#(\xi_L), \tag{6.2}$$

where $c_\gamma > 0$ is a suitable constant depending on $\gamma$. Without loss of generality we may assume that $U = U_1$. Note that for any pair $x, y \in \xi_L$ such that $x \in U$ and $y \in \xi_L \setminus U$ we have $|x - y| = (a_{\gamma L} - x) + \zeta_{\gamma L} + (y - b_{\gamma L})$, the sum of three positive terms. Using $(\sum_{i=1}^3 w_i)^\alpha \geq d(\alpha) \sum_{i=1}^3 w_i^\alpha$, for any $w_i \geq 0$, with $d(\alpha) = 1 \wedge 3^{\alpha-1}$, we can estimate

$$
\begin{aligned}
\Phi_L(\xi) &\leq \frac{\sum_{x \in U} \sum_{y \in \xi_L \setminus U} e^{-|x-y|^\alpha}}{\#(U)} \\
&\leq \frac{1}{c_\gamma L} e^{-d(\alpha)\zeta_{\gamma L}^\alpha} \left( \sum_{x \in U} e^{-d(\alpha)|x-a_{\gamma L}|^\alpha} \right) \left( \sum_{y \in \xi_L \setminus U} e^{-d(\alpha)|y-b_{\gamma L}|^\alpha} \right).
\end{aligned}
\tag{6.3}
$$

Next, we show that each sum in the right-hand side of (6.3) is bounded by $C \log L$ for a suitable constant $C$.

By the assumption $\mathbb{P}(Z_i = 0) = 0$, we can find $\delta > 0$ such that $\mathbb{P}(Z_i \leq \delta) \leq \frac{1}{2}$. We partition $[-\frac{L}{2}, \frac{L}{2}]$ by means of $\delta$-intervals $\Delta_i$ of the form $[a_i, a_{i+1}]$, $a_{i+1} = a_i + \delta$. For any interval $\Delta_i$ we write $N_i = \#(\xi \cap \Delta_i)$ for the number of points in $\Delta_i$. We note that for any integer $K \geq 1$, for every $i$:

$$\mathbb{P}(N_i \geq K) \leq 2^{-K+1}. \tag{6.4}$$

Indeed, we may assume that $\Delta_i = [a, a+\delta]$ for some $a > 0$ (this represents no real loss of generality since the argument when $a \leq 0$ is very similar). Then, we call $X_a = \min\{x \in \xi : x \geq a\}$. Let $\mu_a = \mathbb{P} \circ X_a^{-1}$ denote the law of $X_a$. Then

$$\mathbb{P}(N_i \geq K) = \int_a^{a+\delta} \mu_a(dt) \mathbb{P}(N_i \geq K | X_a = t).$$

From the renewal property we have that for any $t \in \Delta_i = [a, a+\delta]$

$$\mathbb{P}(N_i \geq K | X_a = t) \leq \mathbb{P}(Z_\ell \leq \delta, 1 \leq \ell \leq K-1) \leq 2^{-K+1}.$$

This proves (6.4).

If we take $K = C \log L$ in (6.4) with $C$ suitably large, then a union bound with the Borel–Cantelli lemma shows that we can assume that there is no interval $\Delta_i$ such that $N_i \geq C \log L$. In this case, we can estimate

$$
\begin{aligned}
\sum_{x \in U} e^{-d(\alpha)|x - a_{\gamma L}|^{\alpha}} &\leq \sup_{a \in [-L/2, L/2]} \sum_{x \in \xi_L} e^{-d(\alpha)|x - a|^{\alpha}} \\
&\leq C \log L \sum_{j = -\infty}^{+\infty} e^{-d(\alpha)|j\delta|^{\alpha}} \leq C_2 \log L
\end{aligned}
$$

(6.5)

for a suitable constant $C_2$. The same estimate can be used to treat the second sum in the right-hand side of (6.3). This shows that we can find a constant $C_3$ such that

$$
\Phi_L(\xi) \leq \frac{C_3 (\log L)^2}{L} e^{-d(\alpha)\zeta_{\gamma L}^{\alpha}}.
$$

Finally, it is clear that if $\mathbb{P}(Z_i > \varepsilon) = 1$ for some $\varepsilon > 0$ the argument given above is not needed. Indeed, in this case both sums in (6.3) are finite and therefore (1.23) holds without the $(\log L)^2$ factor in the right-hand side. This ends the proof of Proposition 1.3.

6.2. *General method to estimate* $\Phi_L(\xi)$.  Thanks to Proposition 1.3, the a.s. behavior of $\Phi_L(\xi)$ can be derived from the a.s. behavior of $\zeta_L$. To this aim we recall some results concerning the a.s. asymptotic behavior of maxima of i.i.d. random variables (see [15], Section 4, [11], Section 3.5). Recall that $\psi(t) = \mathbb{P}(Z_1 > t)$.

LEMMA 6.1.  *Define*

$$
M_n := \max(Z_1, Z_2, \ldots, Z_n).
$$

*The following holds:*

(i) *Let* $u_n$ *be a nondecreasing sequence. Then*

$$
\mathbb{P}(M_n \geq u_n, \ i.o.) = 0 \quad or \quad 1
$$

*accordingly as*

$$
\sum_{n=1}^{\infty} \psi(u_n) < \infty \quad or \quad \sum_{n=1}^{\infty} \psi(u_n) = +\infty.
$$

(ii) *Let* $u_n$ *be a nondecreasing sequence such that*

$$
\lim_{n \uparrow \infty} \psi(u_n) = 0,
$$

$$
\lim_{n \uparrow \infty} n\psi(u_n) = \infty.
$$

*Then*

$$\mathbb{P}(M_n \leq u_n, \ i.o.) = 0 \quad or \quad 1$$

*accordingly as*

$$\sum_{n=1}^{\infty} \psi(u_n) e^{-n\psi(u_n)} < \infty \quad or \quad \sum_{n=1}^{\infty} \psi(u_n) e^{-n\psi(u_n)} = +\infty.$$

(iii) *Suppose that* $\mathbb{P}(Z_i \leq x) < 1$ *for all* $x \in \mathbb{R}$ *and set* $\alpha_n := \inf\{y : \psi(y) \leq 1/n\}$. *Then*

$$\mathbb{P}\Big(\lim_{n \uparrow \infty} M_n/\alpha_n = 1\Big) = 1,$$

*if and only if for arbitrary* $k > 1$,

(6.6) $$\sum_{n=1}^{\infty} \psi(k\alpha_n) < \infty.$$

(iv) *Suppose that* $\psi$ *is continuous and that* $\mathbb{P}(Z_i \leq x) < 1$ *for all* $x \in \mathbb{R}$. *Let* $\beta_n$ *be a nondecreasing sequence of positive numbers such that*

$$\mathbb{P}\Big(\lim_{n \uparrow \infty} M_n/\beta_n = 1\Big) = 1.$$

*Then* $\lim_{n \uparrow \infty} \beta_n/\alpha_n = 1$.

As shown in detail in the next subsection, the a.s. asymptotic behavior of $\zeta_L$ is similar to the a.s. asymptotic behavior of $M_n$. Since it is the exponential of $\zeta_L$ that enters in the lower and upper bound of $\Phi_L(\xi)$ in (1.23), it is necessary to know good sequences $u_n, v_n$ such that $u_n \leq M_n \leq v_n$ $\mathbb{P}$-a.s. for $n$ large enough. Thanks to Lemma 6.1(iv), (v), in many cases the right sequences $u_n, v_n$ can be easily guessed, since they must be corrections of $\alpha_n$ (see [11], Chapter 3, for examples and further discussions). These observations are at the heart of Theorem 1.4, which covers several interesting cases already discussed in the Introduction.

6.3. *Proof of Theorem 1.4.* Let us define

(6.7) $n_+(L) := \max\{n : x_n \leq L/2\}, \qquad n_-(L) := \max\{n : x_{-n} \geq -L/2\}.$

Also, for integers $k \geq 0, \ell \geq 1$ we set

(6.8) $$M_{k,\ell} := \max\{Z_{-k}, Z_{-k+1}, \ldots, Z_{\ell-2}, Z_{\ell-1}\}.$$

With these definitions we have

$$\zeta_L = M_{n_-(L), n_+(L)}.$$

We may define $M_{s,t}$ for every $s, t > 0$, by taking the integer parts $M_{s,t} := M_{[s],[t]}$.

From the ergodicity it follows that

$$\lim_{L \to \infty} \frac{n_{\pm}(L)}{L} = \frac{1}{2\mu}, \qquad \mathbb{P}\text{-a.s.,} \tag{6.9}$$

where, as usual $\mu = \mathbb{E}[Z_1]$. For any $\gamma < 1$ and setting $n_L := \frac{L}{\mu}$ we then have

$$[\gamma n_L/2] \le n_{\pm}(L) \le [\gamma^{-1} n_L/2], \qquad \mathbb{P}\text{-a.s.}$$

Therefore

$$M_{1/2\gamma n_L, 1/2\gamma n_L} \le \zeta_L \le M_{1/2\gamma^{-1} n_L, 1/2\gamma^{-1} n_L}, \qquad \mathbb{P}\text{-a.s.} \tag{6.10}$$

Since $n_L$ is deterministic, one can define a deterministic bijective map $\tau : \mathbb{N}_+ \to \mathbb{Z}$ such that, setting $M_n^{(\tau)} = \max(Z_{\tau(1)}, Z_{\tau(2)}, \ldots, Z_{\tau(n)})$, it holds

$$M_{1/2\gamma n_L, 1/2\gamma n_L} = M_{2[\gamma n_L/2]}^{(\tau)}, \tag{6.11}$$

$$M_{1/2\gamma^{-1} n_L, 1/2\gamma^{-1} n_L} = M_{2[\gamma^{-1} n_L/2]}^{(\tau)}. \tag{6.12}$$

The a.s. limiting behavior of $M_n^{\tau}$ can be determined by means of Lemma 6.1, replacing $\{Z_i\}_{i \ge 1}$ with $\{Z_{\tau(i)}\}_{i \ge 1}$. In particular, if $a(\cdot)$ is a function as in Theorem 1.4, satisfying (1.24) it follows from Lemma 6.1(i) that

$$M_{1/2\gamma^{-1} n_L, 1/2\gamma^{-1} n_L} \le a(\gamma^{-1} n_L), \qquad \mathbb{P}\text{-a.s.}$$

Similarly if $b(\cdot)$ is a function as in Theorem 1.4, satisfying (1.25) then Lemma 6.1(ii) implies that

$$M_{1/2\gamma n_L, 1/2\gamma n_L} > b(2[\gamma n_L/2]) \ge b(\gamma^2 n_L), \qquad \mathbb{P}\text{-a.s.}$$

In conclusion, assuming both (1.24) and (1.25) we may estimate, for any $\gamma \in (0, 1)$:

$$b(\gamma^2 n_L) \le \zeta_L \le a(\gamma^{-1} n_L), \qquad \mathbb{P}\text{-a.s.} \tag{6.13}$$

Since $\gamma$ is arbitrarily close to 1 we see that (1.26) is an immediate consequence of Proposition 1.3 and (6.13). This ends the proof of Theorem 1.4.

6.4. *Proof of Corollary 1.5.* Thanks to the simple estimate $\text{gap}(L) \le 2\Phi_L$, the upper bound (1.30) in Corollary 1.5 is an immediate consequence of the upper estimate of Theorem 1.4. Moreover, as discussed in [6] an upper bound of order $L^{-2}$ on $\text{gap}(L)$ is not hard to obtain in our setting. Therefore, we only have to prove the lower bound in (1.29). Let $-\frac{L}{2} \le x_{-n_-(L)} < \cdots <$ $x_{n_+(L)} \leq \frac{L}{2}$ denote the points of $\xi_L$; see (6.7). Let $f$ be an arbitrary vector. The numerator in (1.27) is written as

$$\sum_{x,y \in \xi_L} (f(x) - f(y))^2 = A_L(f) := \sum_{i=-n_-(L)}^{n_+(L)} \sum_{j=-n_-(L)}^{n_+(L)} (f(x_i) - f(x_j))^2.$$

Note that $x_{i+1} - x_i = Z_i$. The denominator in (1.27) therefore satisfies

$$\sum_{x,y \in \xi_L} e^{-|x-y|^\alpha} (f(x) - f(y))^2 \geq B_L(f) := \sum_{i=-n_-(L)}^{n_+(L)-1} e^{-Z_i^\alpha} (f(x_i) - f(x_{i+1}))^2.$$

For any $i < j$ we estimate with the Schwarz inequality

$$(f(x_i) - f(x_j))^2 \leq \sum_{k=i}^{j-1} e^{Z_k^\alpha} \sum_{\ell=i}^{j-1} e^{-Z_\ell^\alpha} (f(x_\ell) - f(x_{\ell+1}))^2.$$

Since $\lim_{n \to \infty} \frac{1}{n} \sum_{k=1}^{n} e^{Z_k^\alpha} = \mathbb{E}[e^{Z_1^\alpha}] < \infty$, $\mathbb{P}$-a.s., using (6.9) we see that for some constant $C_1$ we have, for any $i < j$,

$$(6.14) \qquad (f(x_i) - f(x_j))^2 \leq C_1 L B_L(f), \qquad \mathbb{P}\text{-a.s.}$$

Summing over $i$ and $j$ in (6.14) and using again (6.9) we see that for some constant $C_2$

$$(6.15) \qquad A_L(f) \leq C_2 L^3 B_L(f).$$

From (6.1) it follows that for some constant $C : \gamma(L) \leq CL^2$, $\mathbb{P}$-a.s. This finishes the proof of Corollary 1.5.

## APPENDIX: A SIMPLE ESTIMATE ON RENEWAL POINT PROCESSES

LEMMA A.1. *Suppose that $\xi$ is a renewal point process containing the origin. Then, given $a < b$, there exists positive constants $c, c'$ such that*

$$(A.1) \qquad \mathbb{P}(\xi(a,b) \geq k) \leq ce^{-c'k}, \qquad k \in \mathbb{N}.$$

*Moreover,*

$$(A.2) \qquad \mathbb{E}(\lambda_0(\omega)^k) < \infty, \qquad k \in \mathbb{N},$$

*where $\lambda_x$ is defined by (1.10).*

PROOF. Let us prove (A.1) in the case $a = 0 < b$, the general case can be treated similarly. We fix $u > 0$ such that

$$\gamma := \mathbb{P}(Z_i \leq u) < 1$$

and set $n$ equal to the integer part $[b/u]$. If $Z_1 + Z_2 + \cdots + Z_k < b$ then at least $k - n$ of the variables $Z_1, Z_2, \ldots, Z_k$ are not larger than $u$. Therefore

$$
\mathbb{P}(\xi(0, b) \geq k) = \mathbb{P}(Z_1 + Z_2 + \cdots + Z_k < b) \tag{A.3}
$$

$$
\leq \sum_{j=k-n}^{k} \binom{k}{j} \gamma^j (1 - \gamma)^{k-j}.
$$

The last member is a sum of $n + 1$ addenda. Since $0 \leq \gamma < 1$, for all $j : k - n \leq j \leq k$ it holds that $\gamma^j \leq \gamma^{k-n}$ and $(1 - \gamma)^{k-j} \leq 1$. Moreover

$$
\binom{k}{j} = \binom{k}{k-j} = \frac{k(k-1)(k-2)\cdots(j+1)}{(k-j)!} \leq k^{k-j} \leq k^n.
$$

Hence, for a suitable constant $c > 0$ depending only on $n$,

$$
\sum_{j=k-n}^{k} \binom{k}{j} \gamma^j (1 - \gamma)^{k-j} \leq (n+1) k^n \gamma^{k-n} \leq c \gamma^{k/2} \qquad \forall k \geq 1.
$$

This concludes the proof of (A.1).

The last estimate (A.2) follows immediately from (A.1) and from the trivial bound $\lambda_0(\xi) \leq C \sum_{x \in \mathbb{Z}} \xi([x, x+1)) e^{-c|x|^\alpha}$.   $\square$

Dipartimento di Matematica
Università Roma Tre
Largo S. Murialdo 1
00146 Roma
Italy
E-mail: caputo@mat.uniroma3.it

Dipartimento di Matematica
"G. Castelnuovo"
Università "La Sapienza"
P. le Aldo Moro 2
00185 Roma
Italy
E-mail: faggiona@mat.uniroma1.it